\documentclass[95pt]{article}

\usepackage{verbatim}
\usepackage{hyperref}
\usepackage{amsmath}
\usepackage{amsfonts}
\usepackage{amsthm}
\newtheorem{rem1}{Remark}[section]
\newtheorem{lem1}{Lemma}[section]
\newtheorem{cor1}{Corollary}[section]

\newtheorem{prop1}{Proposition}[section]
\newtheorem{thm1}{Theorem}[section]

\begin{document}
\title{Asymptotics of spectral quantities of Schr\"odinger operators}
\author{T. Kappeler\footnote{Supported in part by the Swiss National Science Foundation}, 
B. Schaad\footnote{Supported in part by the Swiss National Science Foundation}, P. Topalov\footnote{Supported in part by NSF DMS-0901443}}
\maketitle
\begin{abstract} \noindent 
In this paper we provide new asymptotic estimates of the Floquet exponents of Schr\"odinger operators on the circle. 
By the same techniques, known asymptotic estimates of various others spectral quantities are improved.
\end{abstract}    
\section{Introduction}
In this paper we prove asymptotics of various spectral quantities of Schr\"odinger operators $L(q):=-d_x^2+q$ in one space dimension with a potential
 $q$ in $L^2_0(\mathbb{T},\mathbb{R})= \{q\in L^2(\mathbb{T},\mathbb{R});\,\int_0^1q(x) dx=0\}$ 
where $\mathbb{T}$ denotes the one dimensional torus $\mathbb{T}= \mathbb{R}/\mathbb{Z}$. 
The periodic/antiperiodic spectrum of $L(q)$ is real and discrete. 
When listed in increasing order and with their multiplicities the eigenvalues satisfy
\[\lambda_0 <\lambda_1\leq \lambda_2< \lambda_3 \leq \lambda_4< \dots\; .\]
Similarly, the Dirichlet and Neumann spectra of $L(q)$, considered on $[0,1]$, are both real and consist of simple eigenvalues.
We also list them in increasing order 
$\mu_0< \mu_1< \mu_2 <\mu_3< \dots  $ and $\eta_0< \eta_1< \eta_2 <\eta_3< \dots\;$.
Furthermore denote by $M(x,\lambda)$ the fundamental solution of $L(q)$, i.e., the $2\times 2$ matrix valued function satisfying $L(q)M=\lambda M$ and  $M(0,\lambda)= Id_{2\times 2}$, 
\[M(x,\lambda):=\begin{pmatrix} y_1(x,\lambda) & y_2(x,\lambda) \\ y_1'(x,\lambda) & y_2'(x,\lambda)
\end{pmatrix}.\]
When evaluated at  $\lambda=\mu_n$, the Floquet matrix $M(1,\lambda)$ is lower triangular, hence its eigenvalues are given by $y_1(1,\mu_n)$ and $y_2'(1,\mu_n)$. 
By the Wronskian identity, they satisfy $y_1(1,\mu_n)y_2'(1,\mu_n)=1$. By deforming $q$ to the zero  potential along the straight line $tq,\,0\leq t\leq 1,$ one sees that $(-1)^ny_1(1,\mu_n)>0.$ Hence the (normalized) Floquet exponents are given by $\pm \kappa_n$ where
\begin{align}\label{44bis.1}
\kappa_n:= \log \left((-1)^ny_2'(1,\mu_n)\right)=-\log \left((-1)^ny_1(1,\mu_n)\right)
\end{align}
and $\log$ denotes the principal branch of the logarithm.
The $\kappa_n$'s actually play an important role in the spectral theory of Schr\"odinger operators and have applications for the study of KdV  as an integrable PDE -- see e.g. \cite{FM}, \cite{PT}. It turns out that when complemented with the $\mu_n$'s they form a system of canonical coordinates for $L^2_0$. For this reason the $\kappa_n$'s are also referred to as quasi-momenta. The first main result concerns the asymptotics of the
$\kappa_n$'s as $n\rightarrow \infty$. To state it introduce the Sobolev spaces $H^N_0$, $N\in \mathbb{Z}_{\geq 0}$
\[H^N_0 := \left\{p(x)= \sum_{n\neq 0}\hat p_n e^{2\pi inx}\;| \,\left\|p\right\|_N < \infty; \, \hat p_{-n}= \overline{\hat p_{n}}\;\forall n\in \mathbb{Z} \right\}\] with\[\|p\|_N:= \Bigl(\sum_{n\neq 0} |n|^{2N}|\hat p_{n}|^2 \Bigr)^{\frac{1}{2}}.\] 
Note that $H^0_0=L^2_0$.
Furthermore, denote by $H^N_{0,\mathbb{C}}= H^N_0\otimes \mathbb{C} $ the complexification of  $H^N_0$.
For $q\in H^N_{0,\mathbb{C}}$, $L(q)$ is no longer  symmetric with respect to the standard inner product in $L^2_\mathbb{C}\equiv L^2(\mathbb{T},\mathbb{C})$,
\[\langle f,g\rangle = \int_0^1 f(x) \overline{g(x)} dx.\]
The periodic/antiperiodic spectrum as well as the Dirichlet and Neumann spectra are still discrete,
the eigenvalues $(\lambda_n)_{n\geq 0}$, $(\mu_n)_{n\geq 1}$ and $(\eta_n)_{n\geq 0}$ however might be complex valued. We list them with their algebraic multiplicities and in lexicographic ordering, defined for complex numbers $a,b$ by

\[a\preceq b\;  \text{iff}\; [\operatorname{Re}a< \operatorname{Re}b]\; \text{or}\; [\operatorname{Re}a=\operatorname{Re}b\;\text{and}\,\operatorname{Im}a\leq\operatorname{Im}b].\]
 It turns out that there exists a complex neighbourhood $W$ of $L^2_0=  H^0_0$ in $H^0_{0,\mathbb{C}}$ so that on $W$ the Dirichlet eigenvalues $\mu_1,\mu_2,\dots$ are simple and real analytic and   
\[-\log \left((-1)^ny_1(1,\mu_n)\right),\; n=1,2,\dots,\]
define real analytic functions, denoted again by $\kappa_n$. In the sequel, it is convenient to write $\ell^2_n$ for the $n$'th component of a sequence in $\ell^2_\mathbb{C}= \ell^2(\mathbb{N},\mathbb{C})$.
\begin{thm1}\label{4bis.1} Let $N\geq 0$. Then for any $q$ in $W\cap H^N_{0,\mathbb{C}}$, \[\kappa_n= \frac{1}{2\pi n} \left(\langle q, \sin 2\pi n x\rangle+ \frac{1}{n^{N+1}}\ell^2_n\right)\]
uniformly on bounded subsets of $W\cap H^N_{0,\mathbb{C}}$.
\end{thm1}
\begin{rem1} For $N=0$, the asymptotics of Theorem \ref{4bis.1} can be found in \cite{PT} p 60.
 \end{rem1}
\noindent For $q\in H^N_0$, one could expect that the $\kappa_n$'s have an expansion of the form
\[\kappa_n=\sum_{k=0}^{K}\frac{c_k}{n^k}+ \frac{1}{n^{K+1}}\ell^2_n .\]
Surprisingly, Theorem \ref{4bis.1} says that the $c_k$  vanish for any $0\leq k\leq K$ where $K=N$.

\noindent To prove Theorem \ref{4bis.1} we need to define and study special solutions of $L(q)y=\nu^2 y$ for $\nu\in \mathbb{C}$ with $|\nu|$ sufficiently large, obtained by a WKB ansatz,  chosen in such a way that various error terms can be easily estimated -- see Section \ref{sectionspectialsolutions} for details.
It turns out that by the same approach one can improve on the asymptotics of periodic/antiperiodic eigenvalues of $L(q)$ known in the literature.

\begin{thm1}\label{specthm2} Let $q$ be in $H^N_{0,\mathbb{C}}$ with $N\in \mathbb{Z}_{\geq 0}$. Then  \begin{align}\label{11a}\{\lambda_{2n},\lambda_{2n-1}\}= \{m_n \pm\sqrt{\langle q, e^{2\pi inx}\rangle\langle q, e^{-2\pi inx} \rangle +\frac{1}{n^{2N+1}}\ell^2_n } +\frac{1}{n^{N+1}}\ell^2_n\} \end{align} uniformly on bounded subsets of potentials in $H^N_{0,\mathbb{C}}$.
The quantity $m_n$ is of the form 
\begin{align}\label{12a} m_n= n^2\pi^2 +\sum_{2\leq 2j\leq N+1 }c_{2j}\frac{1}{n^{2j}}\end{align}
with coefficients $c_{2j}$ which are independent of $n$ and $N$ and given by integrals of polynomials in $q$ and its derivatives up to order $2j-2$.
\end{thm1}
\begin{rem1} The asymptotic estimate \eqref{11a}, but not of the uniform boundedness of the error
 terms stated above, are due to Marchenko \cite{Ma}.\end{rem1}

\noindent Unfortunately, the asymptotics of Theorem \ref{specthm2} do not lead to satisfactory asymptotic estimates of the sequence $(\tau_n)_{n\geq 1}$ where $\tau_n= (\lambda_{2n}+\lambda_{2n-1})/2$.
In Theorem \ref{reallambdathm} of Appendix B, the asymptotic estimates of Theorem \ref{specthm2} are improved on $H^N_{0}$ and lead to the desired estimates for $\tau_n$. However in the complex case, different arguments have to be used to obtain stronger asymptotics of $\tau_n$.
\begin{thm1}\label{n3.2ter} (i) For any $q\in H^N_0,$ $N\in \mathbb{Z}_{\geq 0},$ 
\begin{align}\label{n17quinto} \tau_n(q)= m_n +\frac{1}{n^{N+1}}\ell^2_n
\end{align} where $m_n$ is given by \eqref{12a} and the error term is uniformly  bounded on bounded sets of potentials in $H^N_0$.

(ii) For any $N\in \mathbb{Z}_{\geq 0}$, there exists an open neighbourhood $W_N\subseteq  H^N_{0,\mathbb{C}}$ of $H^N_0$ so that \eqref{n17quinto} holds on $W_N$ with a locally uniformly bounded error term. 
\end{thm1}
\begin{rem1} We expect that the asymptotics \eqref{n17quinto} hold on all of $H^N_{0,\mathbb{C}}$, and that the error term in \eqref{n17quinto} is bounded on bounded sets of $H^N_{0,\mathbb{C}}$. However, for  the applications in \cite{KST3}, the result as stated suffices.
\end{rem1}
\noindent By the same approach we also obtain a short and self-contained proof of the following asymptotics of the Neumann and of the Dirichlet eigenvalues.

\begin{thm1}\label{specthm} Let $q$ be in $H^N_{0,\mathbb{C}}$ with $N\in \mathbb{Z}_{\geq 0}$. Then 
\begin{align} \label{n15bis}
\eta_n =m_n +\langle q,\cos2\pi nx\rangle+\frac{1}{n^{N+1}}\ell^2_n\\
\label{16a}\mu_n= m_n -\langle q,\cos2\pi nx\rangle +\frac{1}{n^{N+1}}\ell^2_n \end{align} uniformly on bounded subsets of potentials in $H^N_{0,\mathbb{C}}$. 
 Here, $m_n$ is the expression defined in \eqref{12a}.  \end{thm1} 
\begin{rem1}The asymptotics \eqref{16a} of the Dirichlet eigenvalues are due to Marchenko \cite{Ma}. The uniform boundedness 
of the error in \eqref{16a} is shown in \cite{SS}. 
\end{rem1}

\noindent The above theorems are important ingredients in  subsequent work \cite{KST3} on qualitative properties of periodic solutions of KdV  and on the asymptotics of canonically defined normal coordinates, also referred to as Birkhoff coordinates.

\noindent
The paper is organized as follows. In Section \ref{sectionspectialsolutions}, we discuss  special solutions of $L(q)f= \nu^2 f$ which admit an asymptotic expansion as $|\nu| \rightarrow\infty$. Theorem \ref{4bis.1}  is proved in Section \ref{section4bis}, 
Theorem \ref{specthm} in Section \ref{section3}, Theorem \ref{specthm2} in Section \ref{lambdasection} 
and Theorem \ref{n3.2ter} is Section \ref{sectiontau}. In Appendix A we prove results on infinite products needed throughout the paper and in Appendix B we prove improved asymptotics for the periodic/antiperiodic eigenvalues of $L(q)$ for $q$ real valued needed for the proof of Theorem \ref{n3.2ter}.

\section{Special solutions}\label{sectionspectialsolutions}
In this section we prove estimates of special solutions of 
$-y''+qy=\lambda y$ as $|\lambda|\rightarrow \infty$ for potentials $q$ in $H^N_{0,\mathbb{C}}$ with $N\geq 0$ which are needed to derive the claimed asymptotics of various spectral quantities. These solutions are obtained with a WKB ansatz and are a version, suited for our purposes, of solutions introduced and studied by Marchenko \cite{Ma}, p 50 ff. 

As $L(q)= -\partial_x^2+q$ is a differential operator of second order it is convenient to introduce $\nu$ as a new spectral parameter with $\nu^2$ playing the role of $\lambda$. The special solutions considered are denoted by $z_N(x,\nu)$ and defined for $\nu \neq 0$ by 
\begin{align}\label{A.0}
z_N(x,\nu)= y_1(x,\nu^2) +\alpha_N(0,\nu)y_2(x,\nu^2)
\end{align} where $y_1(x,\lambda)$, $y_2(x,\lambda)$ denote the standard fundamental solutions of $-y''+qy=\lambda y$ and $q$ is assumed to be in $H^N_{0,\mathbb{C}}\equiv H^N_0(\mathbb{T},\mathbb{C})$. The function $\alpha_N(x,\nu)$ is given by 
\begin{align}\label{A.0bis}
\alpha_N(x,\nu)= i\nu + \sum_{k=1}^{N}\frac{s_k(x)}{(2i\nu)^k}
\end{align} where 
\begin{align}\label{84bis}s_1(x)=q(x),\quad s_2(x)= -\partial_x q(x) \end{align} 
and, for $2\leq k\leq N$, $s_{k+1}$ is determined by the recursion relation
\begin{align}\label{A.1}
s_{k+1}(x)=-\partial_x s_k(x) - \sum_{j=1}^{k-1}s_{k-j}(x)s_j(x).
\end{align}
(Note that \eqref{A.1} remains true for $N=1$: In this case, the sum in \eqref{A.1} is not present.)
By an induction argument one sees that for any $q\in H^N_{0,\mathbb{C}}$, and any $0\leq k\leq N$, $s_{k+1}$ is a universal isobaric polynomial homogeneous of degree $1+\frac{k}{2}$. Here the adjective 'isobaric' signifies that $q$ is considered to be of degree 1 and differentiation $\partial_x$ of degree 1/2. Furthermore, for any $1\leq k\leq N$,
\begin{align}\label{A.2}
s_{k+1}(x)=(-1)^k\partial_x^k q(x) + S_{k-2}(x)
\end{align}
where $S_{k-2}(x)$ is a polynomial in $q(x),\,\partial_xq(x), \dots, \partial_x^{k-2}q(x)$ with constant coefficients and where $S_{-1}\equiv 0$. As a consequence, $\alpha_N(\cdot, \nu)\in H^1_\mathbb{C}$ and by the Sobolev embedding theorem, $\alpha_N(x,\nu)$ is continuous in $x$ and hence \eqref{A.0} well-defined. Moreover, as $s_j\in H^{N+1-j}_{\mathbb{C}}$ for any $1\leq  j\leq N+1$, $\sum_{j=1}^N s_{N+1-j}(x)s_j(x)$ is in $H^1_\mathbb{C}$ and one may define $s_{N+2}$ by formula \eqref{A.1} as an element in $H^{-1}_\mathbb{C}$
\begin{align}\label{A.1bis}
s_{N+2}(x)=-\partial_x s_{N+1}- \sum_{j=1}^Ns_{N+1-j}s_j=(-1)^{N+1}\partial_x^{N+1} q(x) + H^1_\mathbb{C}.
\end{align}
Clearly, for any $\nu \in \mathbb{C}\setminus \{0\},\, z_N(x,\nu)$ and $z_N(x,-\nu)$  are solutions of $-y''+qy=\lambda y$ with $\lambda=\nu^2$ which both are $1$ at $x=0$ and are linearly dependent solutions iff
\[\alpha_N(0,\nu)-\alpha_N(0,-\nu)=0,
\]i.e. $\nu$ is a zero of the following  polynomial of degree $N+1$
\begin{align*}p_N(\nu)&= (2i\nu)^N\left(\alpha_N(0,\nu)-\alpha_N(0,-\nu)\right)\\&= (2i\nu)^{N+1} + \sum_{k=1}^N s_k(0)\left(1-(-1)^k\right)(2i\nu)^{N-k}.
\end{align*}
It then follows that there exists $\nu_0 > 0$ so that 
\[\left|p_N(\nu)\right|\geq 1 \quad \forall \left|\nu\right|\geq \nu_0.\]
By the Sobolev embedding theorem, the number $\nu_0$ can be chosen uniformly on bounded sets of potentials in $H^N_{0,\mathbb{C}}$. In particular, for $\left|\nu\right|\geq \nu_0$ one has 
\[y_1(x,\nu^2)= \frac{1}{\alpha_N(0,\nu)-\alpha_N(0,-\nu)}\Big(\alpha_N(0,\nu)z_N(x,-\nu)-\alpha_N(0,-\nu)z_N(x,\nu)\Big)\] and 
\[y_2(x,\nu^2)= \frac{1}{\alpha_N(0,\nu)-\alpha_N(0,-\nu)}\left(z_N(x,\nu)-z_N(x,-\nu)\right).\]
Furthermore note that for any $x\in \mathbb{R}$, $z_N(x,\nu)$ is analytic on 
$\mathbb{C}^*\times H^N_{0,\mathbb{C}}$.
As mentioned above, the solutions $z_N(x,\nu)$ are determined by a WKB ansatz,
\[z_N(x,\nu)= w_N(x,\nu)+ \frac{r_N(x,\nu)}{(2i\nu)^{N+1}}\]
with \begin{align}\label{A.2bis}w_N(x,\nu)= \exp\left(\int_0^x\alpha_N(t,\nu)dt\right).\end{align}
By the considerations above it follows that $w_N(\cdot,\nu)$ is in $H^2_\mathbb{C}[0,1]$. As $y_1(\cdot,\nu^2)$ and $y_2(\cdot,\nu^2)$, and hence $z_N(\cdot, \nu)$, are in $H^{N+2}_\mathbb{C}[0,1]$ one then concludes that $r_N(\cdot,\nu)$ is in $H^2_\mathbb{C}[0,1]$ as well.
To study the asymptotics of 
$r_N(x,\nu)$ as $\left|\nu\right| \rightarrow \infty$, first note that 
\[w_N'= \alpha_N w_N \quad \textrm{and}\quad w_N''= (\alpha_N'+\alpha_N^2)w_N.\]
When substituting the latter expression into $-y''+qy=\lambda y$ one obtains 
\begin{align}\label{A.3} \left(-\alpha_N'-\alpha_N^2+ q-\nu^2\right)w_N+ 
\frac{1}{(2i\nu)^{N+1}}\left(-r_N'' +q r_N-\nu^2r_N\right)=0.
\end{align}
Expanding $\alpha_N^2$ in powers of $\nu^{-1}$ leads to
\[\alpha_N^2= -\nu^2+ \sum_{k=0}^{N-1}\frac{s_{k+1}}{(2i\nu)^k}+ \sum_{k=2}^{2N}\frac{1}{(2i\nu)^k}\sum_{1\leq l\leq k-1,\, 1\leq k-l,\, l\leq N}s_{k-l}s_l.\]
Taking into account the identity $s_1=q$ and the relations \eqref{A.1} one then gets \begin{align*}
-\alpha_N'-\alpha_N^2+ q-\nu^2=& -\frac{s_1'+s_2}{2i\nu} -\sum_{k=2}^{N-1}\left(s_k'+s_{k+1}+ \sum_{1\leq l\leq k-1}s_{k-l}s_l\right) \frac{1}{(2i\nu)^k}
\\&-\left(s_N'+ \sum_{1\leq l \leq N-1} s_{N-l}s_l\right)\frac{1}{(2i\nu)^N} \\&-\sum_{k=1}^N\left(\sum_{k\leq l\leq N}s_{N+k-l}s_l\right)\frac{1}{(2i\nu)^{N+k}}\\=& \frac{s_{N+1}}{(2i\nu)^N}- \sum_{k=1}^{N}\left(\sum_{k\leq l\leq N}s_{N+k-l}s_l\right)\frac{1}{(2i\nu)^{N+k}}.
\end{align*}
Substituting this expression into \eqref{A.3} yields 
\begin{align}\label{A.4} -r_N''+qr_N-\nu^2r_N= -2i\nu f_N(x,\nu)
\end{align} where 
\begin{align}\label{84}f_N(x,\nu):=  s_{N+1}w_N-\sum_{k=1}^{N}\left(\sum_{k\leq l\leq N}s_{N+k-l}s_l\right)\frac{w_N}{(2i\nu)^{k}} .\end{align} Recall that
\begin{align}\label{91bis}s_{N+1}(x)-(-1)^N\partial_x^N q(x)\in H^{2}_{\mathbb{C}},\end{align}
$s_{N+k-l}s_l \in H^{1}_{\mathbb{C}} $ for any $k\leq l\leq N$ and $w_N(\cdot, \nu)\in H^2_\mathbb{C}[0,1]$. Hence $f_N(\cdot,\nu) \in L^2_\mathbb{C}[0,1]$. More precisely,
    \[f_N(\cdot,\nu)= (-1)^N\partial_x^Nq \cdot w_N(\cdot,\nu)+ H^1_\mathbb{C}[0,1]\] uniformly on bounded subsets of $H^N_{0,\mathbb{C}}$ and uniformly for $\nu \in \mathbb{C}$ with $|\nu|\geq 1$ and $| \operatorname{Im} \nu |\leq C$.   Further note that
\begin{align}\label{ok93bis}r_N(0,\nu)= (2i\nu)^{N+1}\Big(z_N(0,\nu)-w_N(0,\nu)\Big)=0 \end{align}and \begin{align}\label{ok93ter}
r_N'(0,\nu)=(2i\nu)^{N+1}\Big(z_N'(0,\nu)-\alpha(0,\nu)w_N(0,\nu)\Big)=0.\end{align}
 
 The estimates for $r_N$ are obtained by using that it satisfies the inhomogeneous Schr\"odinger equation \eqref{A.4}. Given $q\in L^2_{0,\mathbb{C}}$ and $\nu\in \mathbb{C}$ denote by $r(x,\nu)$ the unique solution of the initial value problem
 \begin{align}\label{85} -r''+ qr -\nu^2 r=-2i\nu f(x,\nu)
\\\label{86} r(0,\nu)=0 \quad \textrm{and} \quad r'(0,\nu)=0
 \end{align}
 where the inhomogeneous term on the right hand side of \eqref{85} is assumed to be in $L^2_\mathbb{C}([0,1])$ for any value of $\nu$. Note that by assumption
\begin{align}\label{87} Q(1)=0\quad \textrm{where} \quad Q(x)= \int_0^x q(t)dt\quad(0\leq x\leq 1). \end{align}

The solution $r(x,\nu)$  of \eqref{85}-\eqref{86} satisfies the following standard estimates.
\begin{lem1}\label{lem8.0} Let $q\in L^2_{0,\mathbb{C}}$  and $f(x,\nu) =h(x) e^{ix\nu}$ with $h\in L^2_\mathbb{C}[0,1]$. Then for any  $\nu\in \mathbb{C}\setminus \{0\}$ and $0\leq x \leq 1$, the solution $r(x,\nu)$ of \eqref{85}-\eqref{86} satisfies the estimates 
\begin{align*} \left|r(x,\nu) \right| &\leq \left(\frac{R^2}{|\nu|} + \frac{4R^3}{|\nu|^2}\left(1+\frac{1}{|\nu|}\right)\right)\|h\| \\
 \left|r'(x,\nu) \right| &\leq \left(R^2 + \frac{2\left(1+\|q\| \right)R^3}{|\nu|}\left(1+\frac{1}{|\nu|}\right)\right)\|h\|
\end{align*}
where
\[R\equiv R(\nu,q) :=  \exp\left( |\operatorname{Im} \nu |+ \|q\|\right).\]
\end{lem1}
\begin{proof} By the method of the variation of constants, the solution $r(x,\nu)$ has the following integral representation -- see e.g. \cite{PT}, Theorem 2, p 12
\begin{align}\label{94bis} r(x,\nu)= \int_0^x \left(y_1(t,\nu^2)y_2(x,\nu^2)- y_1(x,\nu^2)y_2(t,\nu^2)\right)f(t,\nu) dt 
\end{align} 
where $y_i=y_i (x,\nu^2 ,q) ,\, i=1,2$, denote the fundamental solutions of $-y'' + qy =\nu^2 y$. 
They satisfy the following estimates on $[0,1]\times \mathbb{C} \times L^2_{0,\mathbb{C}}$ -- see e.g. \cite{PT}, Theorem 3, p 13
\begin{align*} \left|y_1(x,\nu^2, q) -\cos \nu x \right| \leq & \frac{R}{|\nu|};&
\left|y_2(x,\nu^2,q) - \frac{\sin \nu x}{\nu}\right| \leq& \frac{R}{|\nu|^2}
\\ \left|y_1'(x,\nu^2, q) +\nu\sin \nu x \right| \leq & \|q\| R ;&
\left|y_2'(x,\nu^2,q) - \cos \nu x \right| \leq& \frac{\|q\| R}{|\nu|}.
\end{align*}
Formula \eqref{94bis} then leads to the following estimates for $0\leq x \leq 1$, $\nu \in \mathbb{C}\setminus \{0\}$, $q\in L^2_{0,\mathbb{C}}$
\begin{align*}\left| r(x,\nu)-\int_0^x \frac{\sin \nu(x-t)}{\nu} f(t,\nu) dt\right| &\leq  \frac{4R^3}{|\nu|^2}\left(1+\frac{1}{|\nu|} \right)\|h\|
\\ \left| r'(x,\nu)-\int_0^x\cos \nu(x-t)f(t,\nu) dt\right| &\leq  \frac{2\left(1+\|q\| \right)R^3}{|\nu|}\left(1+\frac{1}{|\nu|} \right)\|h\|.
\end{align*}
As \[\left|\int_0^x \frac{\sin \nu (x-t)}{\nu} f(t,\nu) dt\right|\leq 
\frac{1}{|\nu|}\int_0^1 e^{2 |\operatorname{Im} \nu |}|h(t)| dt \leq \frac{1}{|\nu|}R^2 \|h\| \] and
\[\left|\int_0^x \cos \nu (x-t) f(t,\nu) dt\right|\leq 
\int_0^1 e^{2 |\operatorname{Im} \nu |}|h(t)| dt \leq R^2 \|h\| \]
it then follows that
\begin{align*}\left| r(x,\nu)\right| &\leq \left(\frac{R^2}{|\nu|} +\frac{4R^3}{|\nu|^2} \left(1+\frac{1}{|\nu|} \right) \right) \|h\| 
\\\left| r'(x,\nu)\right| &\leq \left(R^2 +\frac{2 \left(1+\|q\| \right)R^3}{|\nu|} \left(1+\frac{1}{|\nu|} \right) \right) \|h\| .
\end{align*}
\end{proof}
To obtain the claimed asymptotics of the periodic and Dirichlet eigenvalues, the estimates of Lemma \ref{lem8.0} have to be refined.

\begin{lem1}\label{lemA.1}Assume that $q\in L^2_{0,\mathbb{C}}$. Then for any sequence $\nu_n=n\pi +\frac{1}{n}\ell^2_n, \;n\geq 1,$ the solution $r(x,\nu)$ of \eqref{85}-\eqref{86} satisfies the following estimates: 
\begin{itemize}\item[(i)] If $f(x,\nu)= h(x) e^{ix\nu}$  with $h\in L^2_\mathbb{C}[0,1]$,
\begin{align*}r(1,\pm \nu_n)=&\,(-1)^n \int_0^1h(x) dx + (-1)^{n+1} \int_0^1 h(x) e^{\pm 2in\pi x} dx\\
&\,\pm \frac{(-1)^{n+1}}{2in\pi}\int_0^1 Q(x)h(x)dx +\frac{1}{n}\ell^2_n
\end{align*} 
and 
\begin{align*}r'(1,\pm \nu_n)=& \,\pm in\pi(-1)^n \int_0^1h(x) dx \pm in\pi (-1)^{n} \int_0^1 h(x) e^{\pm 2in\pi x} dx\\
&\,+ \frac{(-1)^{n+1}}{2}\int_0^1 Q(x)h(x)dx +\ell^2_n .
\end{align*} 
\item[(ii)] If $f(x,\nu)= \frac{h(x,\nu)}{\nu^2} e^{i\nu x}$ and the family $h(\cdot, \nu)$ is bounded in $L^2_\mathbb{C}[0,1]$, then 
\[r(1,\pm \nu_n )= O(\frac{1}{n^2}) \quad \text{and} \quad r'(1,\pm\nu_n)= O\left(\frac{1}{n}\right).\]
\end{itemize}
The estimates in (i) and (ii) are uniform on bounded sets of $q's$, $h's$ and $\ell^2_\mathbb{C}$-sequences $\left(\frac{\nu_n-n\pi}{n}\right)_{n\geq 1}$.
\end{lem1}
\begin{proof} As in the proof of Lemma \ref{lem8.0}, the solution $r$ of \eqref{85}-\eqref{86} is written in the following integral form
\begin{align} \label{88} r(x,\nu)= 2i\nu \int_0^x G(x,t;\nu)f(t,\nu)dt
\end{align}
where, with $\lambda= \nu^2$, \[G(x,t;\nu) = y_1(t,\lambda)y_2(x,\lambda)-y_1(x,\lambda)y_2(t,\lambda).\]
As a consequence \[r'(x,\nu) = 2i\nu \int_0^x\partial_x G(x,t; \nu)f(t,\nu) dt.\]
According to \cite{PT} p 14, the solutions $y_i(x,\lambda),\; i=1,2$, and their derivatives  $y_i'(x,\lambda)$ admit the following expansion
\begin{align*}y_1(x,\lambda)=&\; \cos\nu x + \frac{1}{\nu}\int_0^x \sin \nu(x-t)\cdot \cos \nu t \cdot q(t) dt + O(\nu^{-2})\\
y_2(x,\lambda)= &\; \frac{\sin\nu x}{\nu} + \frac{1}{\nu^2}\int_0^x \sin \nu(x-t)\cdot \sin \nu t \cdot q(t) dt + O(\nu^{-3})\end{align*} and  
\begin{align*}y_1'(x,\lambda)=&\; -\nu \sin\nu x + \int_0^x \cos \nu(x-t)\cdot \cos \nu t \cdot q(t) dt + O(\nu^{-1})\\
y_2'(x,\lambda)= &\; \cos \nu x + \frac{1}{\nu}\int_0^x \cos \nu(x-t)\cdot \sin \nu t \cdot q(t) dt + O(\nu^{-2}).\end{align*}
These estimates are uniform on the strip $|\operatorname{Im}\nu| \leq C$, with $C>0$  arbitrary.
 Hence 
\begin{align*}G(x,t;\nu) =&\; \frac{\sin \nu(x-t)}{\nu}+ \frac{1}{\nu^2} \int_0^x \sin \nu(x-s) \cdot \sin \nu(s-t)\cdot q(s) ds 
\\&- \frac{1}{\nu^2} \int_0^t \sin \nu(s-t) \cdot \sin \nu(x-s)\cdot q(s) ds+   O(\nu^{-3}) 
\end{align*} and
\begin{align*}\partial_xG(x,t;\nu) =&\; \cos \nu(x-t)+ \frac{1}{\nu} \int_0^x \cos \nu(x-s) \cdot \sin\nu(s-t)\cdot q(s) ds 
\\&- \frac{1}{\nu} \int_0^t \sin \nu(s-t) \cdot \cos \nu(x-s)\cdot q(s) ds+   O(\nu^{-2}) .
\end{align*}

When substituted into \eqref{88} one gets,
up to an error term which is uniform on bounded sets of $q's$ and $f's$, 
\[r(x,\nu)= I + II +III + O(\nu^{-2}) \] where 
\begin{align*} I\equiv  I(x,\nu) &= 2i \int_0^x \sin \nu(x-t) \cdot f(t,\nu) dt\\
II\equiv  II(x,\nu) &= \frac{2i}{\nu} \int_0^xdtf(t,\nu)\cdot\int_0^x ds \sin \nu(x-s)\cdot \sin \nu(s-t) \cdot  q(s)
\\III\equiv  III(x,\nu) &=- \frac{2i}{\nu} \int_0^xdtf(t,\nu)\cdot \int_0^t ds \sin \nu(x-s)\cdot \sin \nu(s-t) \cdot  q(s)
\end{align*}
After regrouping the terms $\partial_x II$ and $\partial_xIII$, the  derivative $r'(x,\nu)$ can be written in the form
\[r'(x,\nu)= I_1 + II_1 + III_1 + O(\nu^{-1})\]
where 
\begin{align*} I_1 \equiv  I_1(x,\nu)  = \partial_x I(x,\nu) &= 2i\nu \int_0^x \cos \nu(x-t) \cdot f(t,\nu) dt\\
II_1\equiv  II_1(x,\nu)  &= 2i\int_0^xdtf(t,\nu)\cdot\int_0^x ds \cos \nu(x-s)\cdot \sin \nu(s-t) \cdot  q(s)
\\III_1\equiv  III_1(x,\nu)   &=-2i  \int_0^xdtf(t,\nu)\cdot \int_0^t ds \cos \nu(x-s)\cdot \sin \nu(s-t) \cdot  q(s).
\end{align*}
\noindent To prove item (i), each of the three terms are treated seperately. Recall that in (i), $f(x,\nu) $ is of the form $f(t,\nu)= h(t) e^{i\nu t}$. Using that
\[2i \sin \nu(x-t)= e^{i\nu x} e^{-i \nu t}- e^{-i\nu x} e^{i\nu t}\]
term I can be computed to be 
\[I = e^{i\nu x} \int_0^x h(t) dt - e^{-i\nu x}\int_0^x h(t) e^{2i\nu t}dt.\]
As by assumption, $\nu_n= n\pi +\frac{1}{n}\ell_n^2$, one has 
\[e^{\pm i\nu_n t}= e^{\pm in \pi t}\left(1+\frac{1}{n}\ell_n^2\right).\] Thus
\[I(1,\pm \nu_n)=(-1)^n\int_0^1h(t) dt- (-1)^n\int_0^1h(t) e^{\pm 2in\pi t}dt +\frac{1}{n} \ell_n^2.\]
Similarly one shows that 
\[I_1(x,\nu)=i\nu e^{i\nu x}\int _0^x h(t) dt + i\nu e^{-i\nu x}\int_0^x h(t) e^{2i\nu t} dt\]
and therefore 
\[I_1(1,\pm\nu_n)= \pm in\pi(-1)^n \int _0^1 h(t) dt \pm in\pi (-1)^n \int_0^1 h(t) e^{\pm 2in\pi t} dt + \ell^2_n.\]
To treat the term II, write
\begin{align*}2i\sin \nu(x-s)\cdot 2i \sin \nu(s-t)=& e^{i\nu x}e^{-i\nu t}+ e^{-i\nu x}e^{i\nu t}
\\ & -e^{-i\nu x} e^{-i\nu t} e^{2i\nu s}- e^{i\nu x}e^{i\nu t}e^{-2i\nu s}
\end{align*}
to get, with the notation of \eqref{87}, 
\begin{align*}
II(x,\nu) =&\frac{e^{i\nu x}}{2i \nu}\int_0^x h(t) dt\cdot Q(x) + \frac{e^{-i\nu x}}{2i\nu}\int_0^x h(t)e^{2i\nu t} dt\cdot Q(x)\\ &
-\frac{e^{-i\nu x}}{2i\nu }\int_0^x h(t) dt \int_0^x q(s) e^{2i\nu s} ds 
\\ & -\frac{e^{ i\nu x}}{2i\nu }\int_0^x h(t)e^{2i\nu t} dt \int_0^x q(s) e^{-2i\nu s} ds .
\end{align*}
As by assumption \eqref{87}, $Q(1)=0$ one gets, arguing as above
\[ II(1,\pm \nu_n)=\frac{1}{n} \ell_n^2.\]
Similarly, 
\begin{align*}
II_1(x,\nu) =&\frac{e^{i\nu x}}{2}\int_0^x h(t) dt\cdot Q(x) - \frac{e^{-i\nu x}}{2}\int_0^x h(t)e^{2i\nu t} dt\cdot Q(x)\\ &
+\frac{e^{-i\nu x}}{2 }\int_0^x h(t) dt \int_0^x q(s) e^{2i\nu s} ds 
\\ & -\frac{e^{ i\nu x}}{2}\int_0^x h(t)e^{2i\nu t} dt \int_0^x q(s) e^{-2i\nu s} ds 
\end{align*}
leading to 
\begin{align*}
II_1(1,\pm\nu_n) =&\frac{(-1)^n}{2}\int_0^1 h(t) dt\cdot \int_0^1 q(s) e^{\pm2in\pi s}ds \\ &+ \frac{(-1)^{n+1}}{2}\int_0^1 h(t)e^{2in \pi t} dt\cdot  \int_0^1 q(s) e^{\mp 2in \pi s} ds +\frac{1}{n}\ell_n^2.
\end{align*}
The term III is treated similarly, to get
\begin{align*}
III(1,\nu) =&-\frac{e^{i\nu }}{2i \nu}\int_0^1 Q(t) h(t) dt - \frac{e^{-i\nu }}{2i\nu}\int_0^1 Q(t) h(t)e^{2i\nu t} dt\\ &
+\frac{e^{i\nu }}{2i\nu }\int_0^1\int_0^1 1_{[0,t]}(s) q(s) h(t) e^{2i(\nu t-\nu s)} dt ds 
\\ & +\frac{e^{ -i\nu }}{2i\nu }\int_0^1 q(s)\left(\int_s^1 h(t)dt \right)e^{2i\nu s} ds. 
\end{align*}
As $Q(t)h(t)$ and $q(s)\cdot \int_s^1 h(t) dt$ are in $L^2_\mathbb{C}[0,1]$ and $1_{[0,t]}(s) q(s) h(t)$ is in $L^2_\mathbb{C}[0,1]^2$ it then follows that
\[III(1,\pm \nu_n)= \pm \frac{(-1)^{n+1}}{2i n\pi}\int_0^1 Q(t)h(t) dt +\frac{1}{n} \ell_n^2.\] 
Finally, in the same way, one obtains
\begin{align*}
III_1(1,\pm\nu_n) =&\frac{(-1)^{n+1}}{2}\int_0^1 h(t) Q(t) dt +\frac{(-1)^{n+1}}{2} \int_0^1 h(t) \left(\int_0^t q(s) e^{\pm 2i\pi ns}ds\right)dt \\ &+ \frac{(-1)^{n}}{2}\int_0^1 h(t)e^{\pm 2in \pi t} \left(  \int_0^t q(s) e^{\mp 2in \pi s} ds\right)dt \\ & + \frac{(-1)^{n}}{2}\int_0^1 h(t) e^{\pm 2in\pi t}Q(t) dt+\frac{1}{n}\ell_n^2.
\end{align*}

\noindent Towards item (ii) recall that in this case $f(x,\nu)= \frac{h(x,\nu)}{\nu^2}e^{i\nu x}$
and in a straightforward way each of the terms I,II, and III can be bounded pointwise by $O(\nu^{-2})$ whereas each of the terms $I_1,II_1,III_1$ can be bounded pointwise by $O(\nu^{-1})$, leading to the claimed estimates. Going through the various steps of the proof one verifies in a straightforward way that the claimed uniformity of the estimates (i) and (ii) hold.
\end{proof} 

Lemma \ref{lemA.1} will now be applied to get the desired estimates for $r_N(1,\pm \sqrt{\mu_n})$ as $n \rightarrow \infty$. Here $r_N(1,\pm \sqrt{\mu_n})$ is given by \eqref{A.4}-\eqref{ok93ter}. Actually we formulate our results in a slightly more general form. For $q$ in $H^N_{0,\mathbb{C}}$ and $1\leq k\leq N+2$, let us introduce $a_k =\int_0^1 s_k(x) dx$. Recall that $s_{N+2}= (-1)^{N+1}\partial_x^{N+1} q + L^2_\mathbb{C}$. Hence $s_{N+2}$ is in $H_\mathbb{C}^{-1}$ and the integral $\int_0^1 s_{N+2}(x)dx$ is well-defined. By \eqref{A.1bis}
\begin{align}\label{ok98bis} a_{N+2} =-\sum_{j=1}^N \int_0^1 s_{N+1-j}(x)s_j(x) dx.
\end{align}
\begin{prop1}\label{lemA.2}
Let $q$ be in $H^N_{0,\mathbb{C}}$ with $N\geq 0$. Then for any sequence 
$\nu_n=n\pi + \frac{1}{n}\ell_n^2$ one has 
\begin{align*}(i)&  \quad r_N(1,\pm \nu_n)&=&(-1)^na_{N+1} + (-1)^{n+1}(\pm 2in\pi)^N\int_0^1q(x) e^{\pm 2in\pi x}dx 
\\&  &  & \pm \frac{(-1)^{n}}{2in\pi}a_{N+2} + \frac{1}{n}\ell^2_n 
\\ (ii)&
  \quad r_N'(1,\pm \nu_n)&=&\pm in\pi (-1)^na_{N+1} \pm in\pi (-1)^{n}(\pm 2in\pi)^N\int_0^1q(x) e^{\pm 2in\pi x}dx 
\\& & & + \frac{(-1)^{n}}{2}a_{N+2} + \ell^2_n ,
\end{align*}
uniformly for $q$'s in bounded subsets of $H^N_{0,\mathbb{C}}$ and $(\nu_n)_{n\geq 1}$ in sets of sequences such that $(n(\nu_n-n\pi))_{n\geq 1}$ is uniformly bounded in $\ell^2_\mathbb{C}$.
\end{prop1}
\begin{proof} 
(i) Let us first treat the case $N=0$. Then the right hand side of \eqref{A.4} is given by $f_0(t,\nu_n)= q(x)e^{i\nu_n x}$. Hence one has by
 Lemma \ref{lemA.1} (i) with $h=q$, \[r_0(1,\pm\nu_n)= (-1)^{n+1}\int_0^1 q(x) e^{\pm 2in \pi x} dx +\frac{1}{n}\ell_n^2\] where we used that, by assumption, $Q(1)=0$ and that \[\int_0^1 Q(x)q(x)dx= \int_0^1 \frac{1}{2} \partial_x Q(x)^2 dx =0.\]
To apply Lemma \ref{lemA.1} for $N\geq 1$ we write $f_N$ on the right hand side of \eqref{A.4} in the form
\[f_N(x,\nu)= h_1(x) e^{i\nu x}+ \frac{h_2(x)}{2i\nu}e^{i\nu x}+ \frac{h_3(x,\nu)}{\nu^2}e^{i\nu x}.\]
To determine $h_1, \, h_2$, and $h_3$ we need to analyse the expression \eqref{84} defining $f_N(x,\nu)$ in more detail. Recall that by \eqref{A.2bis}, $w_N(x,\nu)= \exp \left(\int_0^x\alpha_N(t,\nu) dt \right)$ and
\[\alpha_N(t,\nu)= i\nu + \frac{q(t)}{2i\nu}+ O(\nu^{-2}).\]
Hence, with $Q(x)=\int_0^x q(t) dt$, 
\[w_N(x,\nu)= e^{i\nu x}\left(1 + \frac{Q(x)}{2i\nu}+ O(\nu^{-2})\right)\]
and
\begin{align}\nonumber f_N(x,\nu) = & s_{N+1}(x)e^{i\nu x} +\left(Q(x)s_{N+1}(x)- \sum_{l=1}^N s_{N+1-l}(x)s_l(x)\right)\frac{e^{i\nu x}}{2i\nu} \\\label{ok98ter} &+ \frac{h_3(x,\nu)}{\nu^2} e^{i\nu x}
\end{align} where $h_3(x,\nu)$  can be explicitly computed from \eqref{84} and $h_3(\cdot,\nu_n)$ is bounded in $L^2_\mathbb{C}[0,1]$ uniformly. 
 Thus by Lemma \ref{lemA.1} one gets
\begin{align*} r_N(1,\pm \nu_n) =& (-1)^n\int_0^1 s_{N+1} dx +(-1)^{n+1}\int_0^1 s_{N+1} e^{\pm 2in\pi x} dx \\& \pm \frac{(-1)^{n+1}}{2in\pi}\cdot \int_0^1 Q\cdot s_{N+1} dx + \frac{1}{n}\ell^2_n
\\& \pm \frac{(-1)^n}{2in\pi}\left(\int_0^1 Q\cdot s_{N+1} dx - \int_0^1 \sum_{l=1}^N s_{N+1-l}\cdot s_l dx\right)\\&
\pm \frac{(-1)^{n+1}}{2in\pi}\int_0^1 Q\cdot s_{N+1}e^{\pm 2in\pi x}dx\\&  \pm \frac{(-1)^n}{2in\pi}\int_0^1\left(\sum_{l=1}^Ns_{N+1-l}\cdot s_l\right)e^{\pm 2in\pi x}dx + \frac{1}{n}\ell^2_n.
\end{align*}
To continue, note that by \eqref{A.2}
\[s_{N+1}(x)= (-1)^N \partial_x^N q(x) + H^2_\mathbb{C}.\]
As $q$ is  periodic, integrating by parts yields
\begin{align*}\int_0^1 s_{N+1}(x) e^{\pm 2in\pi x} dx =& (-1)^N \int_0^1\partial_x^N q(x) e^{\pm 2in\pi x }dx + O(\frac{1}{n^2}) \\ =& (\pm 2in \pi )^N \int_0^1 q(x) e^{\pm 2in\pi x} dx + O(\frac{1}{n^2}).
\end{align*}
Moreover, as $Q(x) s_{N+1}(x)$ and $\sum_{l=1}^N s_{N+1-l}(x)s_l(x)$ are in $L^2_\mathbb{C}$, their sequences of Fourier coefficients are in $\ell_\mathbb{C}^2$. Taking into account identity \eqref{ok98bis} and that the terms containing $\int_0^1Q(x) s_{N+1}(x) dx$ cancel each other we  then  get
\begin{align*} r_N(1,\pm\nu_n) = & (-1)^n \int_0^1 s_{N+1} dx + (-1)^{n+1}(\pm 2in\pi)^N \int_0^1 q(x) e^{\pm 2in\pi x} dx\\& \pm \frac{(-1)^{n}}{2in\pi} a_{N+2} + \frac{1}{n} \ell^2_n.
\end{align*}

\noindent (ii) Again we treat the case $N=0$ first. Then $f_0(t,\nu_n)= q(x)e^{i\nu_n x}$ and hence by Lemma \ref{lemA.1} (i), 
\[r_0'(1,\pm\nu_n) = \pm in\pi (-1)^n \int_0^1 q(x) e^{\pm 2in\pi x}dx + \ell^2_n\]
where we again used that $\int_0^1 Q(x)q(x) dx =0$. As $a_1=0$ and $a_2=0$, the obtained asymptotics coincide with the claimed ones. If $N\geq 1$, then we again use the representation \eqref{ok98ter}  of $f_N(x,\nu)$ to conclude from Lemma \ref{lemA.1} that
\begin{align*}r'(1,\pm\nu_n) =& \pm in\pi (-1)^n a_{N+1} \pm in\pi (-1)^n \int_0^1 (-1)^N \partial_x^Nq(x) e^{\pm 2in\pi x}dx \\&+\frac{(-1)^n}{2}a_{N+2}+ \ell^2_n
\end{align*}
which leads to the claimed asymptotic estimate. Going through the various steps of the proof one verifies in a straightforward way that the claimed uniformity holds.
\end{proof}
At various occasions we will need the following property of the coefficients $s_k$ for $k$ even.
\begin{lem1}\label{lemA.3} For $q$ in $H^N_{0,\mathbb{C}}$ with $N\geq 0$, and $1\leq k\leq N+2,\; a_k= \int_0^1 s_k(x) dx$ is equal to an integral of a polynomial of $q$ and its derivatives up to order $k-3$. Moreover 
\begin{align}\label{100}a_{2k}=\int_0^1 s_{2k}(x) dx =0 \quad \forall\; 2\leq 2k\leq N+2.\end{align}
\end{lem1}

\begin{proof}The first statement follows from the definition \eqref{84bis}-\eqref{A.1} of $s_k$. Indeed by \eqref{A.2}, for $2\leq k\leq N+1,\; a_k= \int_0^1 s_{k-2}(x) dx.$ For the case $k=N+2$ see \eqref{ok98bis}. 
The identities \eqref{100} hold in the case $N=0$ or $N=1$ as $a_2=\int_0^1 -\partial_x q dx=0$ by the definition \eqref{A.1} of $s_2$. It therefore suffices to consider the case where $N\geq 2$. By approximating $q\in H^N_{0,\mathbb{C}}$ with a sequence in $H^{4N}_{0,\mathbb{C}}$, it suffices to proof \eqref{100}  for $q$ in $H^{4N}_{0,\mathbb{C}}$. For $\nu \neq 0$, denote by $Y_{4N}(x,\nu)$ the solution matrix 
\[Y_{4N}(x,\nu)= \begin{pmatrix} z_{4N}(x,-\nu)&z_{4N}(x,\nu)\\z'_{4N}(x,-\nu) & z'_{4N}(x,\nu)
\end{pmatrix}.
\]
As $z_{4N}(0,\nu)=1$ and $z_{4N}'(0,\nu)=\alpha_{4N}(0,\nu)$ one has 
\[\det Y_{4N}(0,\nu) = \alpha_{4N}(0,\nu)-\alpha_{4N}(0,-\nu)= 2i\nu + \sum_{1\leq 2l+1 \leq 4N}\frac{2s_{2l+1}(0)}{(2i\nu)^{2l+1}}.\]
Therefore, $\det Y_{4N}(0,\nu)\neq 0$ for $|\nu|$ sufficiently large.
Furthermore, by the Wronskian identity, $\det Y_{4N}(1,\nu)= \det Y_{4N}(0,\nu)$ and hence,
for $|\nu|$ sufficiently large,
\begin{align}\label{102} \frac{z_{4N}(1,-\nu)z'_{4N}(1,\nu)-z_{4N}(1,\nu)z_{4N}'(1,-\nu)}{\alpha_{4N}(0,\nu)-\alpha_{4N}(0,-\nu)}=1
\end{align}
Write $z_{4N}(1,\pm \nu)$ as a product with an error term as follows. Introduce 
\begin{align*} A:=& \exp\left(i\nu + \sum_{1\leq 2l+1 \leq 4N}\frac{a_{2l+1}}{(2i\nu)^{2l+1}}\right)\\
B:= & \exp \left( \sum_{2\leq 2l \leq 4N}\frac{a_{2l}}{(2i\nu)^{2l}}\right)
\\ R_{\pm} := & \frac{r_{4N}(1,\pm \nu)}{(\pm 2i\nu)^{4N+1}} \quad \text{and}\quad
R_{\pm}' :=  \frac{r_{4N}'(1,\pm \nu)}{(\pm 2i\nu)^{4N+1}}
\end{align*}
As $\alpha_{4N}(x,\nu)$ is periodic in $x$, $\alpha_{4N}(1,\nu)=\alpha_{4N}(0,\nu)$, one concludes 
\begin{align*}z_{4N}(1,\nu)=& A\cdot B+R_+,& z_{4N}(1,-\nu) =& A^{-1} B+ R_-
\\z_{4N}'(1,\nu)=&\alpha_{4N}(0,\nu) A\cdot B+R_+',& z_{4N}'(1,-\nu) =& \alpha_{4N}(0,-\nu)A^{-1} B+ R_-'.
\end{align*}
Substituting these expressions into \eqref{102} yields
\begin{align}\label{104} B^2= 1- \frac{R}{\alpha_{4N}(0,\nu)-\alpha_{4N}(0,-\nu)}
\end{align} where
\begin{align*}R =&\; \alpha_{4N}(0,\nu)R_-AB + R_+'A^{-1}B + R_-R_+' \\
 & -\alpha_{4N}(0,-\nu) R_+A^{-1} B- R_-'AB - R_+ R_-'.
\end{align*}
By Lemma \ref{lem8.0}, applied to $f(x,\nu)= -2i\nu f_{4N}(x,\nu)$ with $f_{4N}(x,\nu)$ defined as in \eqref{84}, one concludes that for $\nu$ real,
\[R_\pm = O\left(\frac{1}{\nu^{4N}}\right) \quad \text{and} \quad R_\pm' = O\left(\frac{1}{\nu^{4N-1}}\right).\]
Furthermore, $\alpha_{4N}(0,\nu)-\alpha_{4N}(0,-\nu)= 2i\nu +O\left(\frac{1}{\nu}\right)$. Taking the logarithm of both sides of \eqref{104} then yields, for $\nu \longrightarrow \infty $, 
\begin{align*}2\sum_{2\leq 2l \leq 4N}\frac{a_{2l}}{(2i\nu)^{2l}}=& \log \left(1- \frac{R}{\alpha_{2N}(0,\nu)-\alpha_{2N}(0,-\nu)}\right)\\ =&\, O\left(\frac{1}{\nu^{4N}}\right)
\end{align*}
which implies that \[a_{2l}=0 \quad \forall 2\leq 2l \leq 4N-2.\]
As $N+2\leq 4N-2$ for $N\geq 2$, the claimed result follows.

\end{proof}

\section{Asymptotics of the $\kappa_n$}\label{section4bis}
In this section we prove Theorem \ref{4bis.1}.
\begin{proof}[Proof of Theorem \ref{4bis.1}]  In section \ref{sectionspectialsolutions}, for $q\in H^N_{0,\mathbb{C}}$, we consider solutions of $-y''+q y= \nu^2 y$ of the form 
\[z_N(x,\pm \nu)= y_1(x,\nu^2)+\alpha_N(0,\pm \nu)y_2(x,\nu^2)\]
where \[\alpha_N(x,\pm\nu)= \pm i \nu +\sum_{k=1}^N \frac{s_k(x)}{(\pm 2i\nu)^k}.\]
(Note that for $N=0$, the latter sum is zero.) Hence for $\nu_n= \sqrt[+]{\mu_n}$ one gets
\[z_N(1,\nu_n)= y_1(1,\mu_n)= z_N(1,-\nu_n).\]
In view of \eqref{B.1} it then follows that 
\begin{align}\label{B1bis}
\kappa_n(q)= -\frac{1}{2}\log \left(z_N(1,\nu_n)z_N(1,-\nu_n)\right).
\end{align}In section \ref{sectionspectialsolutions} we show that 
\[z_N(x,\pm \nu)=w_N(x,\pm \nu)+ \frac{r_N(x,\pm \nu)}{(\pm 2i\nu)^{N+1}}\]
where \[w_N(x,\pm \nu)= \exp \left(\int_0^x \alpha_N(t,\pm\nu) dt\right).\]
Hence $z_N(1,\nu_n)z_N(1,-\nu_n)= I + II + III$ where 
\begin{align*} I =& \exp \left( \sum_{k=1}^N (1+ (-1)^k)\int_0^1 s_k(x) dx \cdot (2i\nu_n)^{-k}\right)
\\ II =&\, w_N(1,-\nu_n) \frac{r_N(1,\nu_n)}{(2i\nu_n)^{N+1}} + w_N(1,\nu_n) \frac{r_N(1,-\nu_n)}{(-2i\nu_n)^{N+1}}
\\ III =&\, r_N(1,\nu_n)r_N(1,-\nu_n)(-1)^{N+1} (2i\nu_n)^{-2N-2}.
\end{align*}
The three terms are analyzed separately. 
Let us begin with I. By Lemma \ref{lemA.3}
\begin{align}\label{B2} (1+ (-1)^k)\int_0^1 s_k(x) dx=0\quad \forall 1\leq k\leq N+2.
\end{align} (These are the cancelations alluded to above.) Therefore  
\begin{align}\label{B3} I=1.
\end{align}
Towards II, note that in view of the assumption $\int_0^1q(x) dx=0 $ one has 
$\mu_n=n^2\pi^2+ \ell^2_n$ and thus 
\[\nu_n= n\pi + \frac{1}{n} \ell^2_n= n\pi \left(1+\frac{1}{n^2}\ell^2_n\right).\]
It implies that
\[(2i\nu_n)^{-N-1}=(2in\pi)^{-N-1} + n^{-N-3}\ell^2_n\]
and \[w_N(1,\pm \nu_n)=(-1)^n\left(1+\frac{1}{n}\ell_n^2\right).\]
Furthermore, by Proposition \ref{lemA.2}, $r_N(1,\pm\nu_n)$ is given by
\begin{align*} (-1)^n \Bigg( a_{N+1}- (\pm 2in\pi )^N \int_0^1 q(x) e^{\pm 2in\pi x} dx
\pm \frac{1}{2in\pi}a_{N+2} + \frac{1}{n} \ell^2_n \Bigg).
\end{align*}
Hence 
\begin{align*} II =& (-1)^n (2in\pi)^{-N-1} \left(r_N(1,\nu_n)+ (-1)^{N+1}r_N(1,-\nu_n)\right) + \frac{1}{n^{N+2}}\ell^2_n \\ 
=& (2in\pi)^{-N-1} \Bigg( (1+ (-1)^{N+1})a_{N+1} - 2i\langle q,\sin 2\pi nx \rangle\\ &
+\frac{1+ (-1)^{N+2}}{2in\pi}a_{N+2} \Bigg) + \frac{1}{n^{N+2}}\ell^2_n.
\end{align*}
Hence again by \eqref{B2} (Lemma \ref{lemA.3}) 
\[(1+ (-1)^{N+2})a_{N+2}=0\quad \text{and} \quad (1+ (-1)^{N+1})a_{N+1}=0.\]
Therefore
\begin{align}\label{B4} II = - \frac{1}{n\pi} \langle q,\sin 2\pi nx \rangle + \frac{1}{n^{N+2}}\ell^2_n.
\end{align}
Finally, 
\begin{align}\label{B5} III= \frac{1}{n^{N+2}}\ell^2_n. \end{align} (If  $N\geq 1$ one has the stronger estimate $III= O(n^{-N-3})$.) Combining \eqref{B3}-\eqref{B5} then yields 
\begin{align*}z_N(1,\nu_n)z_N(1,-\nu_n)= 1- \frac{1}{n\pi} \langle q,\sin 2\pi nx \rangle + \frac{1}{n^{N+2}}\ell^2_n
\end{align*} which, in view of \eqref{B1bis},  leads to 
\begin{align*}\kappa_n =& -\frac{1}{2}\log \left(z_N(1,\nu_n)z_N(1,-\nu_n)\right)\\=&  \frac{1}{2n\pi}\langle q,\sin 2\pi nx \rangle+ \frac{1}{n^{N+2}}\ell^2_n .
\end{align*}
Going through the various steps of the proof one verifies in a straightforward way that the claimed uniformity holds.
\end{proof}

\section{Asymptotics of $\mu_n$ and $\eta_n$}\label{section3}
In this section we prove the asymptotic estimates for the Dirichlet and Neumann eigenvalues as stated in Theorem \ref{specthm} in the introduction. 
\begin{proof}[Proof of Theorem \ref{specthm}]Let us first prove the asymptotics of the Dirichlet eigenvalues.
The main ingredient of the proof are the special solutions $z_N(x,\pm \nu)$ of $-y'' +qy=\nu^2 y$ for $q$ in $H^N_{0,\mathbb{C}}$, constructed in section \ref{sectionspectialsolutions},
\begin{align}\label{b16bis}
 z_N(x,\nu)= \exp\left(\int_0^x\alpha_N(t,\nu) dt\right) + \frac{r_N(x,\nu)}{(2i\nu)^{N+1}}\end{align}
 where $\alpha_N(t,\nu) = i\nu + \sum_{k=1}^N \frac{s_k(t)}{(2i\nu)^k}$ and the functions $s_k(t)$ are given by \eqref{84bis}-\eqref{A.1}. Note that $z_N(0,\pm\nu)=1$ and recall that for $\left|\nu \right|$ sufficiently large, $z_N(x,\nu)$ and $z_N(x,-\nu)$ are linearly independent. Hence $z_N(x,\nu)-z_N(x,-\nu)$ is a scalar multiple of $y_2(x,\nu^2).$ The n'th Dirichlet eigenvalue $\mu_n$ therefore satisfies 

\begin{align}\label{b16ter}
 z_N(1,\nu_n)-z_N(1,-\nu_n) =0 \end{align} where $\nu_n=\sqrt[+]{\mu_n}.$ To analyse \eqref{b16ter} note that by Lemma \ref{lemA.3}, $\exp\left(\int_0^1\alpha_N(t,\nu) dt\right)$ equals 
 \begin{align*}A& = \exp\left(i\nu_n + \sum_{1\leq 2l+1 \leq N} \frac{a_{2l+1}}{(2i\nu_n)^{2l+1}}\right)\\ &=  (-1)^n\exp\left(i(\nu_n-n\pi) + \sum_{1\leq 2l+1 \leq N} \frac{a_{2l+1}}{(2i\nu_n)^{2l+1}}\right)
\end{align*} 
where $a_k=\int_0^1 s_k(t)dt$.  Combining \eqref{b16bis} and \eqref{b16ter} we therefore get the equation
\[A-A^{-1}=R \quad \text{where} \quad R= \frac{r_N(1,-\nu_n)}{(-2i\nu_n)^{N+1}} -\frac{r_N(1,\nu_n)}{(2i\nu_n)^{N+1}}\]
i.e., $A$ satisfies the quadratic equation 
\begin{align}\label{b16*}
 A^2-RA -1 =0. \end{align} 
 As $\mu_n = n^2\pi^2+ \ell^2_n$ one has 
 $\nu_n =n\pi + \frac{1}{n} \ell_n^2$ and hence $A$ is the solution of \eqref{b16*} given by 
 \begin{align}\label{b16quator}
 (-1)^n A= (-1)^nR/2 + \sqrt{1+ R^2/4} = 1+ R^2/8+ O(R^4)  + (-1)^nR/2.\end{align}
According to Proposition \ref{lemA.2} and in view of the asymptotics $\nu_n= n\pi +\frac{1}{n} \ell^2_n$ one has 
\begin{align*}\;-(2i\nu_n)^{N+1} R =&\; r_N(1,\nu_n) + (-1)^{N} r_N(1,-\nu_n)
\\ =&\; (-1)^n (1+(-1)^{N})a_{N+1} + 2 (-1)^{n+1}(2in\pi)^N 
\langle q, \cos 2n\pi x\rangle \\ & + \frac{(-1)^n}{2in\pi}(1+(-1)^{N+1}) a_{N+2} + \frac{1}{n}\ell^2_n.
\end{align*}
As by Lemma \ref{lemA.3}, $(1+(-1)^{N}) a_{N+1}= 2a_{N+1}$ as well as $(1+(-1)^{N+1}) a_{N+2}= 2a_{N+2}$, and 
\begin{align}\label{b16quinto}
 (2i\nu_n)^{-(N+1)}= (2in\pi)^{-(N+1)}(1+\frac{1}{n^2}\ell_n^2) \end{align}
one gets 
\begin{align*}\frac{(-1)^n}{2}R =& -\frac{1}{(2in\pi)^{N+1}}\Big(a_{N+1} - (2in\pi)^N \langle q, \cos 2n\pi x\rangle
 +\frac{1}{2in\pi}a_{N+2} + \frac{1}{n} \ell_n^2
\Big).
\end{align*}
As $a_1=\int_0^1 q(x) dx= 0$ we then conclude that in the case $N=0$,
\[(-1)^n R = \frac{1}{in\pi}\left(\langle q, \cos 2n\pi x\rangle +\frac{1}{n} \ell_n^2 \right)\] and $R^2 =\frac{1}{n^2} \ell_n^2$ whereas for 
$N\geq 1,$ using that $(2in\pi)^N \langle q, \cos 2n\pi x\rangle = \ell_n^2$
\[ R^2 =\frac{1}{(2in\pi)^{N+1}}O\left(\frac{1}{n^2}\right).\] 
Substituting these estimates into \eqref{b16quator} therefore yields in both cases, $N=0$ and $N\geq 1$,
\begin{align*} 
\exp \left( i(\nu_n -n\pi) + \sum_{1\leq 2l+1 \leq N} \frac{a_{2l+1}}{(2i\nu_n)^{2l+1}}\right) =  1+ \frac{(-1)^n}{2}R +O(R^2)
\\ =   1- \frac{1}{(2in\pi)^{N+1}}\Big(a_{N+1} - (2in\pi)^N \langle q, \cos 2n\pi x\rangle
 +\frac{1}{2in\pi}a_{N+2} + \frac{1}{n} \ell_n^2
\Big).
\end{align*}
Taking the principal branch of the logarithm of both sides of the latter identity and multiplying by $-i$ we thus obtain in view of \eqref{b16quinto}, that $\rho_n:= \nu_n-n\pi$ equals
\begin{align} \label{b16sexto}  \rho_n  = & \sum_{1\leq 2l+1 \leq N+2}(-1)^l \frac{a_{2l+1}}{(2\nu_n)^{2l+1}} -\frac{1}{2n\pi} \langle q, \cos 2n\pi x\rangle  + \frac{1}{n^{N+2}}\ell_n^2.
\end{align}
Unfortunately, $\nu_n$ appears also on the right hand side of \eqref{b16sexto}.
 To address this issue, we follow an approach found by Marchenko \cite{Ma}. Let 
 \begin{align}\label{n112} F(z) = \sum_{1\leq 2l+1 \leq N+2}\frac{(-1)^l}{2^{2l+1}} a_{2l+1}z^{2l+1}\end{align} and write 
 \[\frac{1}{\nu_n}= \frac{1}{n\pi + \rho_n} =  \frac{1/n}{\pi + \rho_n/n}.\]
 We approximate $F(\frac{1}{\nu_n})$ by approximating $\rho_n$ by $\rho(1/n)$ in the above expression where $\rho$ is an analytic function so that near $z=0$, \[\rho(z)- F\left(\frac{z}{\pi + z \rho(z)}\right)=0.\] To find $\rho$ introduce 
 \[G(z,w):= w- F\left(\frac{z}{\pi + z w}\right).\]
Note that $G(0,0)=0$ and $\partial_w G(0,0) =1$. Hence by the implicit function theorem there exists near $z=0$ a unique analytic function $\rho=\rho(z)$ so that $\rho(0)=0$ and $G(z,\rho(z))=0$ for $z$ near $0$.
Note that $F$ is an odd function; hence 
\[G(-z,-w)=-w+F\left(\frac{z}{\pi + z w}\right)= -G(z,w)\]
and as a consequence, $G(-z,-\rho(z))=0$ near $z=0$.
On the other hand, $G(-z,\rho(-z))=0$ and therefore, by the uniqueness of $\rho(z)$, one has $\rho(-z) = -\rho(z)$. It follows that $\rho$ has an expansion of the form
\begin{align}\label{n113}\rho(z) = \sum_{k=0}^\infty b_{2k+1}z^{2k+1}.\end{align}
The coefficients $b_{2k+1}$ can be  computed recursively from the identity 
\[\rho(z) = F\left(\frac{z}{\pi + z \rho(z)}\right).\]
In this way one sees that for any $k\geq 0 $, $b_{2k+1}$ is a polynomial in the coefficients $a_{2l+1}$ of $F$ with $0\leq l\leq k.$
The Taylor expansion of $F(z)$ at $z_n= \frac{1}{n\pi +\rho(\frac{1}{n})}$ with Lagrange's remainder term reads 
\begin{align}\label{n114} F\left(\frac{1}{\nu_n}\right)= \rho\left(\frac{1}{n}\right) + F_n\cdot \left(\rho\left(\frac{1}{n}\right) -\rho_n\right)\end{align}
where 
\begin{align}\label{n115} F_n= \int_0^1 F'\left(z_n+t\left(\frac{1}{\nu_n}-z_n\right)\right)dt \cdot \frac{1}{\nu_n \cdot (n\pi +\rho(1/n))}= O(\frac{1}{n^2})
\end{align} as $\rho(0)=0$ and $\rho_n=\frac{1}{n}\ell^2_n \rightarrow 0$ as $n\rightarrow \infty$.
Substracting $\rho(1/n)$ on both sides of the identity \eqref{b16sexto} then
yields 

 \[(1+F_n)\cdot(\rho_n -\rho(1/n))= -\frac{1}{2n\pi}\langle q, \cos 2n\pi x\rangle +\frac{1}{n^{N+2}}\ell_n^2 \] or
\[\rho_n-\rho(1/n) = -\frac{1}{2n\pi}\langle q, \cos 2n\pi x\rangle +\frac{1}{n^{N+2}}\ell_n^2.\]
Thus we have shown that 
\[\nu_n =n\pi + \sum_{1\leq 2k+1 \leq N+2} b_{2k+1}\frac{1}{n^{2k+1}}-\frac{1}{2n\pi}\langle q, \cos 2n\pi x\rangle +\frac{1}{n^{N+2}}\ell_n^2. \]
By taking squares on both sides of the latter identity we obtain the asymptotics \eqref{16a} with the claimed properties of the expression $m_n$ of \eqref{12a}. Going through the arguments of the proof one verifies the stated uniformity property of the error term in \eqref{16a}.

\noindent The asymptotic estimates for the Neumann  eigenvalues $(\eta_n)_{n\geq 0}$  are derived in a similar way as the ones for the Dirichlet eigenvalues. Note that the special solutions $z_N(x,\pm\nu)$ satisfy $z_{N}(0,\pm\nu)=1$ and $z_N'(0,\pm \nu)=\alpha_N(0,\pm\nu)$.  Hence $y(x,\nu):=\alpha_N(0,-\nu)z_N(x,\nu)-\alpha_N(0,\nu)z_N(x,-\nu)$ satisfies $y'(0,\nu)=0$. As $\alpha_N(0,\pm \nu)= \pm i\nu + \sum_{1\leq k\leq N}\frac{s_k(0)}{ (\pm 2i\nu)k }$, for $|\nu|$ sufficiently large,
\[y(0,\nu)= \alpha_N(0,-\nu)-\alpha_N(0,\nu)= -2i\nu + O\left(\frac{1}{\nu}\right)\neq 0.\]
Therefore $y(x,\nu)$ is parallel to $y_1(x,\nu^2)$. Again it is convenient to introduce $\nu_n=\sqrt{\eta_n}=n\pi+ \frac{1}{n}\ell^2_n$. The $n$'th Neumann eigenvalue $\eta_n$ is then characterized by
\begin{align}\label{n200} \alpha_N(0,-\nu_n)z'_N(1,\nu_n)-\alpha_N(0,\nu_n)z'_N(1,-\nu_n)=0.
\end{align} 
To analyze \eqref{n200} note that 
\[z'_N(1,\pm\nu_n)= \alpha_N(1,\pm \nu_n)A^{\pm 1}+ \frac{r'_N(1,\nu_n)}{(\pm 2i\nu_n)^{N+1}},
\]where $A^{\pm 1}= 
\exp \left( \pm i\nu_n + \sum_{1\leq 2l+1 \leq N} \frac{a_{2l+1}}{(\pm2i\nu_n)^{2l+1}}\right)$. 
As $\alpha_N(x,\pm \nu_n)$ is 1-periodic in $x$, \eqref{n200} reads 
\begin{align*} &\alpha_N(0,-\nu_n)\alpha_N(0,\nu_n)A- \alpha_N(0,-\nu_n)\alpha_N(0,\nu_n)A^{-1} \\=& \alpha_N(0,\nu_n)\frac{r_N'(1,-\nu_n)}{(-2i\nu_n)^{N+1}}-\alpha_N(0,-\nu_n)\frac{r_N'(1,\nu_n)}{(2i\nu_n)^{N+1}}.
\end{align*}
Thus $A$ satisfies the quadratic equation 
\[A^2-2RA -1=0\]
where \[R=\frac{1}{2\alpha_N(0,-\nu_n)}\frac{r'_N(1,-\nu_n)}{(-2i\nu_n)^{N+1}}-\frac{1}{2\alpha_N(0,\nu_n)} \frac{r_N'(1,\nu_n)}{(2i\nu_n)^{N+1}}.\]
As $\nu_n= n\pi + \frac{1}{n}\ell^2_n$, we write \[(-1)^nA= \exp \left(  i(\nu_n-n\pi) + \sum_{1\leq 2l+1 \leq N} \frac{a_{2l+1}}{(2i\nu_n)^{2l+1}}\right)\] leading to 
\begin{align}\label{n201}(-1)^nA= (-1)^nR + \sqrt[+]{1+R^2}= 1+(-1)^nR + O(R^2).\end{align}
By Proposition \ref{lemA.2}, $(-1)^nr'_N(1,\pm\nu_n)$ equals
\begin{align*} \pm in \pi a_{N+1} \pm in\pi (\pm2in\pi)^N\int_0^1q(x)e^{\pm 2\pi inx}dx + \frac{a_{N+2}}{2}+\ell^2_n.
\end{align*}
Further 
\begin{align*}\alpha_N(0,\pm \nu_n)= \pm i \nu_n + \sum_{1\leq k\leq N}\frac{s_k(0)}{(\pm 2i\nu_n)^{k}}
= \pm in\pi \left(1+O\left(\frac{1}{n^2}\right)\right),
\end{align*}
and thus
\begin{align*} \frac{(-1)^nr'_N(1,\pm\nu_n)}{\alpha_N(0,\pm \nu_n)}=& a_{N+1} + (\pm 2in\pi)^N\int_0^1q(x)e^{\pm 2\pi inx}dx
+ \frac{a_{N+2}}{\pm2in\pi}+ \frac{1}{n}\ell^2_n
\end{align*} yielding the asymptotic estimate
\begin{align*} (-1)^nR =  \frac{1}{2} \frac{(-1)^{N+1}}{(2in\pi)^{N+1}}a_{N+1}-\frac{1}{2} \frac{1}{2in\pi} \int_0^1 q(x)e^{-2\pi inx}dx +\frac{a_{N+2}}{2} \frac{(-1)^{N+2}}{(2in\pi)^{N+2}}
\\  -\frac{1}{2} \frac{1}{(2in\pi)^{N+1}}a_{N+1} -\frac{1}{2} \frac{1}{2in\pi} \int_0^1 q(x)e^{2\pi inx}dx -\frac{a_{N+2}}{2} \frac{1}{(2in\pi)^{N+2}}
+ \frac{1}{n^{N+2}}\ell^2_n.
\end{align*}
Using that by Lemma \ref{lemA.3},
\[a_{N+1}(1+(-1)^N)= 2a_{N+1}\quad \text{and}\quad a_{N+2}(1+(-1)^{N+1})= 2a_{N+2}, \]
one then gets \[(-1)^nR = -a_{N+1}\frac{1}{(2in\pi)^{N+1}}-\frac{1}{2in\pi}\langle q,\cos 2\pi nx \rangle -a_{N+2}\frac{1}{(2in\pi)^{N+2}}+\frac{1}{n^{N+2}}\ell^2_n.\]
In view of $a_1=a_2=0$ one then gets for any $N\geq 0$
\[R^2= \frac{1}{n^{N+2}}\ell^2_n.\]
Substituting the latter two estimates into \eqref{n201} one concludes 
\begin{align*} &\exp \left(  i(\nu_n-n\pi) + \sum_{1\leq 2l+1 \leq N} \frac{a_{2l+1}}{(2i\nu_n)^{2l+1}}\right)\\
 =& 1- a_{N+1}\frac{1}{(2in\pi)^{N+1}}- \frac{1}{2in\pi } \langle q,\cos 2\pi nx \rangle- a_{N+2}\frac{1}{(2in\pi)^{N+2}}+ \frac{1}{n^{N+2}}\ell^2_n
\end{align*}
and taking the principal branch of the logarithm on both sides one gets after multiplying by $-i$ for any $N\geq 0$,
\begin{align*} \rho_n:= \nu_n-n\pi = \sum_{1\leq 2l+1 \leq N+2} \frac{(-1)^la_{2l+1}}{(2\nu_n)^{2l+1}}+\frac{1}{2n\pi } \langle q,\cos 2\pi nx \rangle+ \frac{1}{n^{N+2}}\ell^2_n.
\end{align*}
Arguing as in the case of the Dirichlet eigenvalues one then concludes that 
\begin{align*} \nu_n =& n\pi +\sum_{1\leq 2l+1 \leq N+2} \frac{b_{2l+1}}{n^{2l+1}} 
+ \frac{1}{2n\pi } \langle q,\cos 2\pi nx \rangle+ \frac{1}{n^{N+2}}\ell^2_n.
\end{align*}
Squaring both sides of the latter identity yields the claimed asymptotics for the Neumann eigenvalues. Going through the arguments of the proof one verifies the stated uniformity property of the error term. 
\end{proof}

\section{Asymptotics of $\lambda_n$}\label{lambdasection}
 The main purpose of this section is to prove Theorem \ref{specthm2}.
\begin{proof}[Proof of Theorem \ref{specthm2}] As for the proof  of Theorem \ref{specthm}, the main ingredient are the special solutions $z_N(x,\pm\nu_n)$ of $-y''+qy= \nu^2y$ for  potentials $q\in H^N_{0,\mathbb{C}}$, constructed in section \ref{sectionspectialsolutions},
\[z_N(x,\nu)= \exp \left(\int_0^x\alpha_N(t,\nu)dt\right)+ \frac{r_N(x,\nu)}{(2i\nu )^{N+1}}\]
where $\alpha_N(t,\nu)= i\nu + \sum_{1\leq k \leq N}\frac{s_k(t)}{(2i\nu)^k}.$ (Without further reference, we use the notation  introduced in section \ref{sectionspectialsolutions}.) Recall from section \ref{sectionspectialsolutions} that for $|\nu|$ sufficiently large, $z_N(x,\nu)$ and $z_N(x,-\nu)$ are linearly independent. Denote by $Y_N(x,\nu)$ the solution matrix 
\[Y_N(x,\nu)= \begin{pmatrix} z_N(x,-\nu) & z_N(x,\nu) \\ z'_N(x,-\nu) & z'_N(x,\nu)\end{pmatrix}\]
and recall that $r_N(0,\nu)=0$ and $r'_N(0,\nu)=0$ so that 
\[Y_N(0,\nu)= \begin{pmatrix}1&1 \\ \alpha_N(0,-\nu) & \alpha_N(0,\nu).
\end{pmatrix}\]
The large periodic eigenvalues of $-d_x^2+q$ on the interval $[0,1]$ are thus given by the zeros of the characteristic function \[\chi_p(\nu):= \det \left(Y_N(1,\nu)-Y_N(0,\nu)\right)\]
whereas the large antiperiodic eigenvalues of $-d_x^2+q$ on the interval $[0,1]$ are given by the zeroes of \[\chi_{ap}(\nu):= \det\left(Y_N(1,\nu)+Y_N(0,\nu)\right).\]
The two cases are treated in a similar fashion and hence we concentrate on the periodic case only. Recall that for $q=0$, they are given by $\lambda_0=0,\; \lambda_{4n}= \lambda_{4n-1}= (2n\pi)^2$. 
For arbitrary $q$ one then knows that the large periodic eigenvalues are $\lambda_{2n},\; \lambda_{2n-1}$ with $n$ large and \textit{even}. It is convenient to introduce the notation 
\[\nu_n^+= \sqrt[+]{\lambda_{2n}},\quad \nu_n^-= \sqrt[+]{\lambda_{2n-1}}\] and to write $\nu_n$ if we do not need to specify our choice among $\nu_n^+$ and $\nu_n^-$. Let us now compute the asymptotics of $\chi_p(\nu)$. First note that $\alpha_N(1,\nu)= \alpha_N(0,\nu)$ as $q$ is $1$-periodic, leading to the formula
\[z'_N(1,\nu)=\alpha_N(0,\nu)\exp\left(\int_0^1 \alpha_N(t,\nu)dt\right)+ \frac{r'_N(1,\nu)}{(2i\nu)^{N+1}}.\]
Then we have
\[z_N(1,\pm\nu_n)=A^{\pm1} +R_{\pm}\] and \[z'_N(1,\pm\nu_n)=\alpha^{\pm}A^{\pm1} +R'_{\pm}\]
where \[A^{\pm 1}:= \exp\left(\pm\int_0^1\alpha_N(t,\nu_n)dt\right)= \exp\left(\int_0^1\alpha_N(t,\pm\nu_n)dt\right)\]
(the latter identity follows from Lemma \ref{lemA.3})
\[R_{\pm}:= \frac{r_N(1,\pm\nu_n)}{(\pm2i\nu_n)^{N+1}}\quad R'_{\pm}:= \frac{r'_N(1,\pm\nu_n)}{(\pm2i\nu_n)^{N+1}}\]
and \[\alpha^{\pm}:= \alpha_N(0,\pm \nu_n)= \pm i\nu_n+ \sum_{k=1}^N\frac{s_k(0)}{(\pm 2i\nu_n)^k}.\]
Hence $\chi_p(\nu_n)$ equals \begin{align*} \left(A^{-1} +R_--1\right)\left(\alpha^+A+R'_+ -\alpha^+\right)
 - \left(\alpha^-A^{-1} + R'_- -\alpha^-\right)\left(A+R_+-1\right)\end{align*}
or
\[\chi_p(\nu_n)= \xi A+ \zeta+ \eta A^{-1}
\]
where 
\begin{align*}\xi =& \alpha^+R_-- R'_--\alpha^+ +\alpha^-\\
\zeta =& 2\alpha^+-2\alpha^-+R'_--\alpha^+R_- - R'_+ + \alpha^-R_+ + R_-R'_+ - R'_-R_+\\
\eta =& -\alpha^-R_+ +R'_+  -\alpha^+ +\alpha^- .
\end{align*}
Note that 
\[\zeta = - (\xi+\eta)+R \quad \text{with} \quad R =R_-R'_+- R'_-R_+.\]
As $\chi_p(\nu_n)=0$, one has $\xi A+ \zeta +\eta A^{-1}=0$ or 
\begin{align}\label{nc200} A= \frac{\xi+\eta -R+\epsilon \sqrt{(\xi+\eta-R)^2-4\xi\eta}}{2\xi}
\end{align}with an appropriate choice of the sign $\epsilon =\epsilon_n^{\pm}\in \{\pm 1\}$. We want to estimate the terms on the right hand side of \eqref{nc200}. 
By Proposition \ref{lemA.2} and the assumption that $n$ is \textit{even}, $(-1)^n=1$ and hence
\[r_N(1,\pm\nu_n)=a_{N+1}-e_n^{\pm} \pm \frac{1}{2in\pi}a_{N+2} + \frac{1}{n} \ell^2_n\]
where \[e_n^{\pm}= (\pm2in\pi)^N\int_0^1 q(x) e^{\pm 2in\pi x}dx.\]
Similarly, the asymptotics of $r'_N(1,\pm\nu_n)$ are given by
\[r'_N(1,\pm\nu_n)= \pm in\pi a_{N+1} \pm in\pi e_n^{\pm}+ \frac{1}{2}a_{N+2} +\ell^2_n.\]
Furthermore, as $a_1=\int_0^1 q dx=0$ by assumption, and $\nu_n= n\pi +\frac{1}{n} \ell^2_n$
\[\alpha^{\pm}= \pm in\pi +\frac{1}{n}\ell^2_n.\]
These asymptotics yield the following estimates
\begin{align}\label{nc201} \xi =& \alpha^- -\alpha^+ + \frac{-a_{N+1}+a_{N+2}/(2in\pi)+\ell^2_n/n}{(-2in\pi)^N} \\\label{nc202}
\eta =&\alpha^- -\alpha^+ + \frac{a_{N+1}+a_{N+2}/(2in\pi)+\ell^2_n/n}{(2in\pi)^N}.
\end{align} 
Furthermore, $R=R_-R'_+-R'_-R_+$ has an expansion of the form
\begin{align}\label{nc203} R =\frac{i}{(2n\pi)^{2N+1}}\left(a_{N+1}^2-e_n^+e_n^-+ \frac{1}{n}\ell^2_n \right).
\end{align} 
Now we are ready to estimate the terms on the right hand side  of \eqref{nc200}. By Lemma \ref{lemA.3}, $a_{N+1}=0 \;[a_{N+2}=0]$ for $N$ odd [even]. Hence for any $N$, $(1+(-1)^{N+1})a_{N+1}=0$ and $(1+(-1)^{N})a_{N+2}=0$. As a consequence, \eqref{nc201}-\eqref{nc202} yield 
\[\xi +\eta= 2(\alpha^- -\alpha^+) +\frac{1}{n^{N+1}}\ell^2_n\] and 
\[\eta -\xi = \frac{1+(-1)^N}{(2in\pi)^N}a_{N+1} + \frac{1+(-1)^{N+1}}{(2in\pi)^{N+1}}a_{N+2}+\frac{1}{n^{N+1}}\ell^2_n.\]
Furthermore, if $N=0$, then use $a_1=\int_0^1 q(x) dx=0$ to conclude that by \eqref{nc201}, $R=\frac{1}{n}\ell^2_n$, whereas for $N=1,$ the fact that $a_2=0$ leads to the estimate $R=\frac{1}{n^2}\ell^2_n.$ For $N\geq 2$, one gets from \eqref{nc201} that $R=\frac{1}{n^{N+1}}\ell^2_n$. Altogether we  have established that for any $N\geq0$
\[R= \frac{1}{n^{N+1}}\ell^2_n\] 
and thus \[\xi +\eta-R= 2(\alpha^--\alpha^+) + \frac{1}{n^{N+1}}\ell^2_n. \]
To estimate the square root in \eqref{nc200}, the term  $R$ will play a role. First note that
\[ (\xi +\eta -R)^2-4\xi \eta = (\eta-\xi)^2 -2(\xi+\eta)R+R^2.\]
As $a_{N+1}\cdot a_{N+2}=0$ for any $N\geq0$ by Lemma \ref{lemA.3}, we get \begin{align*}(\eta-\xi)^2=&\frac{(1+(-1)^N)^2}{(2in\pi)^{2N}}a_{N+1}^2 + \frac{1}{n^{2N+1}}\ell^2_n\\ 
-2(\xi+\eta)R=&\left(4(\alpha^+-\alpha^-)+\frac{1}{n^{N+1}}\ell^2_n\right)\frac{(-1)^{N+1}}{(2in\pi)^{2N+1}}\left(a_{N+1}^2-e_n^+e_n^-+\frac{1}{n}\ell^2_n \right) 
\\ =&\frac{(-1)^{N+1}4}{(2in\pi)^{2N}}\left(a_{N+1}^2-e_n^+e_n^-\right)+\frac{1}{n^{2N+1}}\ell^2_n .\end{align*}
Clearly $R^2 = \frac{1}{n^{2N+1}}\ell^2_n$ and thus 
\begin{align*} (\xi +\eta-R)^2 -4\xi \eta =& \frac{(1+(-1)^N)^2-(-1)^{N}4}{(2in\pi)^{2N}}a_{N+1}^2 +\frac{(-1)^{N}4}{(2in\pi)^{2N}}e_n^+e_n^- \\&+\frac{1}{n^{2N+1}}\ell^2_n.
\end{align*}
Using once more that $a_{N+1}=0$ for $N$ odd it follows that for any $N\geq 0$,
\[\left((1+(-1)^N)^2-(-1)^{N}4\right)a_{N+1}^2=0\]
leading to\[ \sqrt{(\xi+\eta-R)^2 -4\xi\eta}= \frac{2}{(2n\pi)^N}\sqrt{e_n^+e_n^-+\frac{1}{n}\ell^2_n}. \]
Finally we need to estimate $1/2\xi$. As $(\alpha^- -\alpha^+)= -2in\pi +O(\frac{1}{n}) $ it follows from \eqref{nc201} that 
\[2\xi = 2(\alpha^- -\alpha^+) \left(1+\frac{(-1)^Na_{N+1}}{(2in\pi)^{N+1}}-\frac{(-1)^Na_{N+2}}{(2in\pi)^{N+2}}+ \frac{1}{n^{N+2}}\ell^2_n\right).\]
As  $a_{N+1}=0$ for $N=0$, we get from the Taylor expansion $(1+x)^{-1}= 1-x+ O(x^2)$, for any $N\geq 0$,
\[\frac{1}{2\xi}= \frac{1}{2(\alpha^--\alpha^+)}\left(1+ \frac{(-1)^{N+1}a_{N+1}}{(2in\pi)^{N+1}}+\frac{(-1)^Na_{N+2}}{(2in\pi)^{N+2}}+\frac{1}{n^{N+2}}\ell^2_n\right).\]
Combining the estimates obtained so far and substituting them into the identity \eqref{nc200} we get for any $N\geq 0$, after dividing nominator and denominator of the right hand side by $2(\alpha^--\alpha^+)$,
\begin{align*} A=& \left(1\pm \frac{i\epsilon_n}{(2n\pi)^{N+1}}\sqrt{e_n^+e_n^-+ \frac{1}{n}\ell^2_n} +\frac{1}{n^{N+2}}\ell^2_n\right)\\ 
& \cdot \left(1- \frac{(-1)^Na_{N+1}}{(2in\pi)^{N+1}}+\frac{(-1)^Na_{N+2}}{(2in\pi)^{N+2}}+\frac{1}{n^{N+2}}\ell^2_n\right) \\
=& 1+ \frac{(-1)^{N+1}a_{N+1}}{(2in\pi)^{N+1}}+ \frac{(-1)^{N+2}a_{N+2}}{(2in\pi)^{N+2}}\\ 
&\pm \frac{i\epsilon_n}{(2n\pi)^{N+1}}\sqrt{e_n^+e_n^-+ \frac{1}{n}\ell^2_n} + \frac{1}{n^{N+2}}\ell^2_n.
\end{align*}
Here we used that $\left(\sqrt{e_n^+e_n^-+ \frac{1}{n}\ell^2_n}\right)_{n\geq 1}$ is in $\ell^2$. Taking the principal branch of the logarithm of  both sides and taking into account that $\log (1+x)= x+ O(x^2)$ as well as $\nu_n-n\pi= \frac{1}{n}\ell^2_n$ and $e^{in\pi}=1$ as $n$ is \textit{even} one gets 
\begin{align*} i(\nu_n-n\pi)+ \sum_{1\leq 2l+1\leq	N} \frac{a_{2l+1}}{(2i\nu_n)^{2l+1}} =& 
\frac{(-1)^{N+1}a_{N+1}}{(2in\pi)^{N+1}}+ \frac{(-1)^{N+2}a_{N+2}}{(2in\pi)^{N+2}} \\ 
& \pm \frac{i\epsilon_n}{(2n\pi)^{N+1}}\sqrt{e_n^+e_n^-+ \frac{1}{n}\ell^2_n} + \frac{1}{n^{N+2}}\ell^2_n.
\end{align*}
Using once more that by Lemma \ref{lemA.3}, $a_k=0$ for $k$ even we get, after multiplying both sides by $-i$, the following estimate for $\rho_n:= \nu_n-n\pi$ \begin{align*} \rho_n= \sum_{1 \leq 2l+1\leq N+2}(-1)^{l}\frac{a_{2l+1}}{(2\nu_n)^{2l+1}} +\frac{\epsilon_n}{(2n\pi)^{N+1}}\sqrt{e_n^+e_n^-+ \frac{1}{n}\ell^2_n} + \frac{1}{n^{N+2}}\ell^2_n.
\end{align*}
Arguing as in the proof of Theorem \ref{specthm} (section \ref{section3}), one has by \eqref{n112}-\eqref{n115}
\begin{align*}\sum_{1 \leq 2l+1\leq N+2}(-1)^{l}\frac{a_{2l+1}}{(2\nu_n)^{2l+1}}= \rho\left(\frac{1}{n}\right) + F_n\cdot\left(\rho\left(\frac{1}{n}\right)-\rho_n\right)
\end{align*} 
with $\rho(z)$ given by \eqref{n113}  and $F_n$ by \eqref{n115}.
Therefore \[(1+F_n)\cdot \left(\rho_n-\rho\left(\frac{1}{n}\right)\right)=\frac{\epsilon_n}{(2n\pi)^{N+1}}\sqrt{e_n^+e_n^-+ \frac{1}{n}\ell^2_n} + \frac{1}{n^{N+2}}\ell^2_n. \]
By \eqref{n115}, $F_n= O(1/n^2)$ and thus
\[\nu_n=n\pi + \sum_{1\leq 2l+1\leq N+2} b_{2l+1} \frac{1}{n^{2l+1}} +\frac{\epsilon_n}{(2n\pi)^{N+1}}\sqrt{e_n^+e_n^-+ \frac{1}{n}\ell^2_n} + \frac{1}{n^{N+2}}\ell^2_n.\] 
Squaring the latter expression yields
\[\nu_n^2= m_n+\frac{\epsilon_n}{(2n\pi)^{N}}\sqrt{e_n^+e_n^-+ \frac{1}{n}\ell^2_n} + \frac{1}{n^{N+1}}\ell^2_n\]
with $m_n$ given by \eqref{12a}. Going through the arguments of the proof one verifies that the error terms are uniformly bounded on bounded subsets of potentials in $H^N_{0,\mathbb{C}}$.
\end{proof}

 We have also the following result on the 
coefficients $c_{2j}$ in \eqref{12a}. 

\begin{cor1}\label{cor3.0} Let $N\in \mathbb{Z}_{\geq 0}$. Then for any $2\leq 2j \leq N$,  $c_{2j}$ is a spectral invariant on $H^N_{0,\mathbb{C}}$, i.e., for any two potentials $p,q$  in $H^N_{0,\mathbb{C}}$ so that $-d_x^2 +p$ and $-d_x^2 +q$ have the same periodic spectrum, one has $c_{2j}(p)=c_{2j}(q)$. In addition, if $N+1$ is even (otherwise $c_{N+1}$ vanishes on $H^N_{0,\mathbb{C}}$) 
there exists an open neighbourhood $\tilde W_N \subseteq H^N_{0,\mathbb{C}}$ of $H^N_0$ so that  on $\tilde W_N$, $c_{N+1}$ is a spectral invariant as well. 
\end{cor1}

\begin{proof}[Proof of Corollary \ref{cor3.0}] Let $p,q$ in $H^N_{0,\mathbb{C}}$ be isospectral, i.e.,  $-d_x^2+p$ and $-d_x^2+q$ have the same periodic spectrum. By the asymptotic estimates \eqref{11a} of Theorem \ref{specthm2} we have $c_{2j}(p)=c_{2j}(q)$ for any $2\leq 2j\leq N$. The case of the coefficient $c_{N+1},\; N+1$ even, is more subtle as the factor in the asymptotic  estimate \eqref{11a} containing Fourier coefficients of $q$, is of comparable size. 

By \cite{KaP}, Theorem 11.10 and Theorem 11.11,  there exists an open neighbourhood $\tilde W_N \subseteq H^N_{0,\mathbb{C}}$ of $H^N_0$ so that any two isospectral potentials $p,q$ in 
 in $\tilde W_N$  can be approximated by isospectral finite gap potentials. In particular it follows that there exist sequences  
$(p_l)_{l\geq 1},(q_l)_{l\geq 1}$
in $\tilde W_N$ with the following properties
\begin{itemize} 
\item[(i)]  $(p_l)_{l\geq 1},(q_l)_{l\geq 1}\subset H^k_{0,\mathbb{C}}$ for any $k\in  \mathbb{Z}_{\geq 0}$;
\item[(ii)] $\lim_{l\rightarrow \infty}p_l =p$, $\lim_{l\rightarrow \infty}q_l =q$ in $H^N_{0, \mathbb{C}}$; 
\item[(iii)]  the periodic spectra of $-d_x^2 +p_l$ and $-d_x^2 +q_l$ coincide for any $l$ in  $\mathbb{Z}_{\geq 1}$.
\end{itemize}

%and introduce
%\[w:=\Phi(p) = (w_n)_{n\geq 1}= (u_n,v_n)_{n\geq 1}\;\text{and}\;z:= \Phi(q)= (z_n)_{n\geq 1}= (x_n,y_n)_{n\geq 1}.\]
%Denote by $(w^{(l)})_{l\geq 1}$ the approximating sequence defined by
%\[ w_n^{(l)}=w_n\;\text{if}\; n\leq l\quad \text{and}\quad w_n^{(l)}=(0,0)\;\text{if}\; n>l,\] 
% and by $(z^{(l)})_{l\geq 1}$ the corresponding sequence obtained from $z$. 

%\noindent Clearly $\lim_{l\rightarrow \infty}w^{(l)}=w$ and $\lim_{l\rightarrow \infty}z^{(l)}=z$. 
%For $l$ sufficiently large, $l\geq l_0$, $w^{(l)}$ and $z^{(l)}$ are contained in $\Phi(\tilde W_N)$ as $\Phi(\tilde W_N)$ is open in $\mathfrak{h}^{N+1/2}_\mathbb{C}$. 

%For those $l$, $p_l=\Phi^{-1} (w^{(l)})$ and $q_l= \Phi^{-1}(z^{(l)})$ are in $W_N$. By construction, $p_l$ and $q_l$ are isospectral and $\lim_{l\rightarrow \infty}p_l =p$, $\lim_{l\rightarrow \infty}q_l =q$. Furthermore, by the definition \eqref{n82bis} of the Birkhoff map together with \eqref{X4}, Proposition \ref{lemX.10} and Corollary \ref{corX.15} 
%\[(\tau_n-\mu_n)(p_l)=0\quad \forall n>l\quad \text{and}\quad (\tau_n-\mu_n)(q_l)=0\quad \forall n>l.\]
%It then follows from \cite{Poschel}, Theorem 5, (with $\delta_n$ in Theorem 5 given by $\delta_n=\mu_n-\tau_n$) that $p_l,q_l$ are in any Sobolev space $H^s_{0,\mathbb{C}}$.
\noindent By (i) and the first part of this proof it then follows that 
\[c_{N+1}(p_l)=c_{N+1}(q_l) \quad \forall l\geq 1.\]
As $c_{N+1}$ is the integral of a polynomial in $q$ and its derivatives up to order $N-1$, it follows that, $c_{N+1}(p)=\lim_{l\rightarrow\infty} c_{N+1}(p_l)$ and $c_{N+1}(q)=\lim_{l\rightarrow\infty} c_{N+1}(q_l)$. Hence $c_{N+1}(p)=c_{N+1}(q)$ as claimed.
\end{proof}

\begin{rem1} Denote by $\Delta(\lambda)\equiv \Delta(\lambda,q)$ the discriminant of $-d_x^2+q$, i.e., $\Delta(\lambda)= y_1(1,\lambda) +y_2'(1,\lambda)$. Using Corollary \ref{cor3.0} one can prove that the Poisson bracket $\{c_{N+1}, \Delta(\lambda)\}= \int_0^1\partial_qc_{N+1} \partial_x\partial_q\Delta(\lambda)$ vanishes on $H^N_0$ for any $\lambda \in \mathbb{C}$. As $c_{N+1}$ and $\Delta(\lambda)$ are defined on all of $H^N_{0,\mathbb{C}}$ and are analytic there $\{c_{N+1}, \Delta(\lambda)\}=0$ on $H^N_{0,\mathbb{C}}$ for any $\lambda \in \mathbb{C}$.
\end{rem1}

\section{Proof of Theorem \ref{n3.2ter}}\label{sectiontau}
The main purpose  of this section is to prove the asymptotics of the $\tau_n$'s  stated in Theorem \ref{n3.2ter}.
\begin{proof}[Proof of Theorem \ref{n3.2ter} (i)] Theorem \ref{n3.2ter} (i) is a direct consequence of the asymptotics  of the periodic eigenvalues stated in Section \ref{asymptoticsofperiodicreal}, Theorem \ref{reallambdathm}.
\end{proof} 
In the case of complex valued potentials the arguments are more involved.
The strategy  to prove Theorem \ref{n3.2ter} (ii)  is to try to find for any given element $q$ in $H^N_{0,\mathbb{C}}$ an isospectral  potential $p$ with the property that the periodic spectrum consists of the disjoint union of the Dirichlet and Neumann spectrum. Such a $p$ is found with the help of the so called Birkhoff map (cf \cite{KaP} for a detailed construction) for any $q$ in $H^N_{0,\mathbb{C}}$ sufficiently close to the real subspace $H^N_0$. The asymptotics of $\tau_n$ are then obtained by applying Theorem \ref{specthm} and Corollary \ref{cor3.0}. First we need to establish some auxiliary results.
Introduce the following subset of $L^2_{0,\mathbb{C}}$ 
\[E:= \{q\in L^2_{0,\mathbb{C}}|\; spec_D(-d_x^2+q) \subseteq  spec_p(-d_x^2+q)\}.\]

We begin by examining some properties of $E\cap L^2_0.$ It turns out that for $q$ in $E\cap L^2_0$,
each Neumann eigenvalue is a periodic one as well. To prove this fact we first need to establish the following result for even potentials. We say that $q\in L^2_{0,\mathbb{C}}$ is even if $q(x)=q(1-x)$ for a.e. $0<x<1.$ 
\begin{lem1}\label{nD.1} Assume that $q\in L^2_{0,\mathbb{C}}$ is even. Then each Dirichlet and each Neumann eigenvalue of $-d_x^2 +q$ is also a periodic eigenvalue.
\end{lem1}
\begin{proof}
Consider first the Dirichlet eigenvalues $(\mu_n)_{n\geq 1}$. For any $n\geq 1$, the function  $g(x):= y_2(1-x, \mu_n)$ satisfies the equation ($0<x<1$) 
\begin{align*} -g''(x) +q(x) g(x)=& - y_2''(1-x,\mu_n) +q(1-x)y_2(1-x,\mu_n)\\
=& \mu_n y_2(1-x,\mu_n)
\end{align*}
where we used that by assumption $q(x)=q(1-x)$ for a.e. $0<x<1$. As $g(0)= y_2(1,\mu_n)$ it then follows that $g(x)=g'(0)y_2(x,\mu_n) $ for any $0\leq x\leq 1$. But $g'(0)= -y_2'(1,\mu_n)$ and therefore 
\[-1= g'(1) = g'(0) y_2'(x,\mu_n)|_{x=1}= -y'_2(1,\mu_n)^2.\]
As a consequence $y'_2(1,\mu_n)=\pm 1$, implying that $\mu_n$ is a periodic eigenvalue of $-d_x^2+q$ (when considered of [0,2]).
For the Neumann eigenvalues $\eta_n,\;n\geq 0,$ one argues similarly. Recall that $y'_1(1,\eta_n)=0$. Consider $h(x):= y_1(1-x,\eta_n)$. Then $h'(0)=0$ and a.e. $0<x<1$, 
\[-h''(x) +q(x)h(x) =\eta_n h(x).\]
Hence $h(x)= h(0)y_1(x,\eta_n)$, or, when evaluated at $x=1$, $1= y_1(1,\eta_n)^2$, again implying that $\eta_n$ is a periodic eigenvalue.
\end{proof}
Lemma \ref{nD.1} allows us to prove the following result for elements in $E\cap L^2_0$, mentioned above. 
\begin{lem1}\label{nD.2}
For any $q\in E\cap L^2_0$, the Neumann spectrum $(\eta_n)_{n\geq 0}$ of $-d_x^2 + q$ is contained in the periodic spectrum as well and one has 
\[\eta_0=\lambda_0, \quad \{\eta_n,\mu_n\}= \{\lambda_{2n-1},\lambda_{2n}\}\quad \forall n\geq 1.
\]
\end{lem1}
\begin{proof}As $\mu_n(q)$ is assumed to be a periodic eigenvalue, one has $y'_2(1,\mu_n)=(-1)^n$ and hence $\kappa_n=\log (-1)^n y'_2(1,\mu_n)=0$ for any $n\geq 1$. By \cite{PT}, Lemma 3.4, it then follows that $q$ is even and thus by Lemma \ref{nD.1} the Neumann eigenvalues $(\eta_n)_{n\geq 0}$ are also periodic ones. Furthermore, as $q$ is assumed to be real, one has $\eta_0 \leq \lambda_0$ and, for any $n\geq  1,\; \lambda_{2n-1}\leq \mu_n, \eta_n \leq \lambda_{2n}$. Hence $\eta_0=\lambda_0$ and  $\{\eta_n,\mu_n\}= \{\lambda_{2n-1},\lambda_{2n}\}$ for any $n\geq 1$ as claimed.
\end{proof}

 The Birkhoff map is defined on  an open neighbourhood $W$ of $L^2_0$ and takes values in
$\mathfrak{h}^{1/2}_\mathbb{C}$, where  for $\alpha\in \mathbb{R}_{\geq 0}$, 
$\mathfrak{h}^\alpha_\mathbb{C} = \ell^{2, \alpha}_{\mathbb{C}}\times \ell^{2, \alpha}_{\mathbb{C}}$ and
\[\ell^{2, \alpha}_{\mathbb{C}} = 
\left\{u=(u_k)_{k\geq 1}\subseteq \mathbb{C}| \; 
\|u\|_{\ell^{2,\alpha}} = \left(\sum_{k=1}^\infty k^{2\alpha}|u_k|^2\right)^{1/2}<\infty \right\}.\] 
Similarly we let $\mathfrak{h}^\alpha:= \ell^{2, \alpha}\times \ell^{2, \alpha}$ denote  the corresponding spaces of real sequences.
 Denote by $ (x_n(q),y_n(q))_{n\geq 1}$ the image $\Phi(q)$.  First we characterize the image of $E\cap W$ by $\Phi$. For this purpose introduce the closed subspace $Z$ of $\mathfrak{h}^{1/2}_\mathbb{C}$,
\[Z:= \{ (x_k,y_k)_{k\geq 1}\in \mathfrak{h}^{1/2}_\mathbb{C} |\; y_k=0 \; \forall k\geq 1\}.\] 
\begin{lem1}\label{nD.3} Elements of $E\cap W$ are mapped by $\Phi$ to $Z\cap \Phi(W)$ and 
\[\Phi|_{E\cap L^2_0}: E\cap L^2_0 \to Z\cap \mathfrak{h}^{1/2}\] is bijective.
\end{lem1}
\begin{proof}
By a straightforward computation it follows from the definition of the Birkhoff map, analyzed in \cite{KaP},  that 
 $\Phi(E\cap L^2_0)\subset Z\cap \mathfrak{h}^{1/2}$. Moreover, as 
$\Phi|_{L^2_0}: L^2_0 \to \mathfrak{h}^{1/2}$ is a diffeomorphism the restriction of $\Phi$ to $L^2_0\cap E$ is $1-1$. To show that $\Phi(L^2_0\cap E) = \mathfrak{h}^{1/2} \cap Z$, let $(x_k,0)_{k\geq 1}$ be an arbitrary element in $\mathfrak{h}^{1/2}\cap Z,$ and define $p= \Phi^{-1}((x_k,0)_{k\geq 1})$. Denote by $q$ the potential in $L^2_0\cap E$ with the same periodic spectrum as $p$ and determined uniquely by the conditions
\begin{align*}\mu_k(q) = 
\begin{cases} \lambda_{2k} & \text{if} \; x_k\leq 0 \\
\lambda_{2k-1} &\text{if}\; x_k>0.
\end{cases}
\end{align*}
By a straightforward calculation one then concludes that $\Phi(q)=\Phi(p)$ and hence $p=q$.

\end{proof}
According to \cite{KaP}, for any $N\in \mathbb{Z}_{\geq 0}$, there exists an open neighbourhood $\tilde W_N$ of $H^N_0$ in $W\cap H^N_{0,\mathbb{C}}$ such that $ \Phi|_{\tilde W_N}: \tilde W_N \to \Phi(\tilde W_N)\subseteq \mathfrak{h}^{N+1/2}_\mathbb{C}$ is a diffeomorphism. We denote $\Phi(\tilde W_N)$ by $\tilde V_N$.
\begin{prop1}\label{nD.4}For any $N\in \mathbb{Z}_{\geq 0}$ and any $q\in \Phi^{-1}(\tilde V_N \cap Z)$.
\[\eta_0(q)= \lambda_0(q) \quad \{\mu_n(q),\eta_n(q)\}= \{\lambda_{2n}(q), \lambda_{2n-1}(q)\}\; \forall n\geq 1.\]
\end{prop1}
\begin{proof}
As $\Phi|_{\tilde W_N}$ is a real analytic diffeomorphism onto its image $\tilde V_N$ and $Z\cap \tilde V_N$ is a real analytic submanifold of $\tilde V_N$ it follows that $\Phi^{-1}(\tilde V_N \cap Z)$ is a real analytic submanifold of $\tilde W_N$. Both, Dirichlet and Neumann eigenvalues are all simple for $q\in W$. It follows that they are real analytic functions on $W$ and so are their restrictions to $\tilde W_N$. Using that the discriminant $\Delta(\lambda,q)$, given by $\Delta(\lambda,q)= y_1(1,\lambda)+ y_2'(1,\lambda)$, is analytic on $\mathbb{C}\times L^2_{0,\mathbb{C}}$ it then follows that the compositions
\[F_n: \tilde W_N \to \mathbb{C},\; q\mapsto \Delta(\mu_n(q),q)^2-4\quad (n\geq 1)\]
and
\[G_n: \tilde W_N \to \mathbb{C},\; q\mapsto \Delta(\eta_n(q),q)^2-4\quad (n\geq 0)\]
are analytic. One verifies easily that they are real on $\tilde W_N\cap H^N_0.$ By the definition of $E$, $F_n \; (n\geq 1)$ vanishes on $E\cap H^N_0$, and by Lemma \ref{nD.2}, so does $G_n\;(n\geq 0)$.
Furthermore, by Lemma \ref{nD.3}, one concludes that $E\cap H^N_0= \Phi^{-1}(Z\cap \mathfrak{h}^{N+1/2})$.
Hence, being analytic, the functions $F_n\; (n\geq 1)$ and $G_n\; (n\geq 0)$ vanish on  $\Phi^{-1}(\tilde V_N\cap Z)$. By the choice of $W$ (see end of introduction), for any $q\in W$ and $n\geq 1$, the eigenvalues $\mu_n(q), \eta_n(q),\lambda_{2n-1}(q),\lambda_{2n}(q)$ are contained in an isolating neighbourhood
$U_n$. The $U_n$'s are pairwise disjoint and none of them contains $\lambda_0$ or $\eta_0$. This implies the claimed statement.
\end{proof}
Using Proposition \ref{nD.4} we want to find a neighbourhood $W_N \subseteq \tilde W_N$ of $H^N_0$ so that for any $q\in W_N$ there exists an isospectral potential $p\in W_N\cap E$. First we establish the following elementary result. For $\beta\in \mathbb{R}_{\geq 0}$ define 
\[ \ell^{1,\beta}_\mathbb{C} = \left\{v=(v_k)_{k\geq 1}\subseteq \mathbb{C}|
 \; \|v\|_{\ell^{1,\beta}} = \sum_{k=1}^\infty k^{\beta}|v_k|<\infty \right\}.\]
Introduce the map 
\[Q:\ell^{2,\alpha}_\mathbb{C}\to \ell^{1,2\alpha}_\mathbb{C},\quad (u_k)_{k\geq 1}\mapsto
Q((u_k)_{k\geq 1})= (u_k^2)_{k\geq 1}.\]
Obviously one has $\|Q(u)\|_{\ell^{1,2\alpha}}=\|u\|_{\ell^{2,\alpha}}^2$. Denote by $B_{\ell^{1,\beta}}(v;\epsilon)\, [B_{\ell^{2,\alpha}}(u;\epsilon)]$ the open ball of radius $\epsilon$ in
 $\ell^{1,\beta}_{\mathbb{C}}\,[\ell^{2,\alpha}_\mathbb{C}]$, centered at $v\in \ell^{1,\beta}_\mathbb{C}\, [u \in \ell^{2,\alpha}_\mathbb{C}].$
\begin{lem1}\label{nD.60} For any $\alpha \geq 0,\; \epsilon >0$, and $u \in \ell^2_{\alpha,\mathbb{C}}$
\[Q(B_{\ell^{2,\alpha}}(u;\sqrt{\epsilon}))\supseteq  B_{\ell^{1,2\alpha}}(Q(u);\epsilon).\]
\end{lem1}
\begin{proof}
Let $u\in \ell^{2,\alpha}_\mathbb{C}$ be arbitrary and consider $v=(v_k)_{k\geq 1}\in B_{\ell^{1,2\alpha}}(Q(u);\epsilon).$ We would like to find 
$(h_k)_{k\geq 1}\in \ell^{2,\alpha}_\mathbb{C}$ so that for any 
$k\geq 1$, \[(u_k+h_k)^2= v_k\;\text{implying that} \;\|(u_k+h_k)^2-u_k^2\|_{\ell^{1,2\alpha}}< \epsilon.\]
For any $k\geq 1$ we obtain the following quadratic equation  for $h_k$
\[h_k^2+ 2u_kh_k-b_k\quad \text{where}\quad b_k=v_k-u_k^2.\]
Choose the sign $\sigma_k\in \{\pm 1\}$ of the root in $h_k=-u_k + \sigma_k\sqrt{v_k^2+b_k}$ in such a way that \[|h_k|= \min_{\pm}|-u_k\pm \sqrt{u_k^2+b_k}|.\]
Then \[h_k^2 \leq |(-u_k+\sigma_k\sqrt{u_k^2+b_k})(-u_k-\sigma_k\sqrt{u_k^2+b_k})|=|b_k|
\] and therefore\[\|(h_k)_{k\geq 1}\|_{\ell^{2,\alpha}}\leq \|v-Q(u)\|_{\ell^{1,2\alpha}}<\epsilon.\]
\end{proof}
The following lemma will allow us to deal with complex potentials.
\begin{lem1}\label{nD.6}
For any $z^0\in \mathfrak{h}^{N+1/2}$ there exists  a neighbourhood $V(z^0)$ of $z^0$ in $\tilde V_N$ so that for any $z=(x_k,y_k)_{k\geq 1}\in V(z_0)$ there exists an element $(u_k,0)_{k\geq 1}$ in $\tilde V_N$ with $u_k^2=x_k^2+y_k^2$ for any $k\geq 1$.
\end{lem1}
\begin{proof} Let $z^0=(x_k^0,y_k^0)_{k\geq 1}$ be an arbitrary real sequence in $\mathfrak{h}^{N+1/2}$.
Define $u^0=  (u_k^0)_{k\geq 1}$ in $\ell^{2,N+1/2}$ by setting $u^0_k=\sqrt[+]{v_k^0}$ for any $k\geq 1$
where $v_k^0=(x_k^0)^2+(y_k^0)^2$. As $(u_k^0,0)_{k\geq 1}\in \mathfrak{h}^{N+1/2}$ and $\tilde V_N(=\Phi (\tilde W_N))$ is open in $\mathfrak{h}^{N+1/2}_\mathbb{C}$, we can choose $\epsilon>0$ so that 
\[\{ (u_k,0)_{k\geq 1}|\; (u_k)_{k\geq 1}\in B_{\ell^{2,N+1/2}}(u^0;\sqrt{\epsilon})\}\subseteq \tilde V_N.
\] Now consider the map \[P: \tilde V_N \to \ell^{1,2N+1}_\mathbb{C},\quad (x_k,y_k)_{k\geq 1}\mapsto (x_k^2+y_k^2)_{k\geq 1}. \] Clearly, $P$ is continuous and therefore there exists a neighbourhood $V(z^0)$ of $z^0$  in $\tilde V_N$ so that 
\[P(V(z^0))\subseteq B_{\ell^{1,2N+1}}(v^0;\epsilon).\]
By Lemma \ref{nD.60}, 
$ B_{\ell^{1,2\alpha}}(Q(v^0);\epsilon)\subseteq Q(B_{\ell^{2,\alpha}}(v^0;\sqrt{\epsilon}))$
and hence 
the neighbourhood $V(z^0)$ has the claimed property.
\end{proof}
Having made these preparations we are ready to prove Theorem \ref{n3.2ter}.
\begin{proof}[Proof of Theorem \ref{n3.2ter} (ii)]
Let $W_N\subseteq \tilde W_N$ be the open neighbourhood of $H^N_0$ in $W\cap H^N_{0,\mathbb{C}}$ given by
\[W_N = \bigcup_{q\in H^N_0} \Phi^{-1}(V(\Phi(q)))\]
where $V(\Phi(q)) \subseteq \tilde V_N$ is the neighbourhood of Lemma \ref{nD.6}. Then for any $q\in W_N$ there exists an element $p\in \tilde W_N$ so that $\Phi(p)= (u_k,0)_{k\geq 1}\in \tilde V_N$  satisfies
$u_k^2= x_k^2+y_k^2$ where $\Phi(q))=(x_k,y_k)_{k\geq 1}.$ Hence $p$  and $q$ are isospectral. In particular $\tau_n(q)=\tau_n(p)$ for any $n\geq 1$. By Proposition \ref{nD.4} $\{\mu_n(p),\eta_n(p)\}=\{\mu_n(q),\eta_n(q)\}$ for any $n\geq 1$ and hence by Theorem \ref{specthm}
\begin{align*} \tau_n(p)=m_n(p)+\frac{1}{n^{N+1}}\ell^2_n.\end{align*}
As by Corollary \ref{cor3.0}, $m_n(p)=m_n(q)$, one then gets 
\begin{align}\label{n260} \tau_n(q)=m_n(q)+\frac{1}{n^{N+1}}\ell^2_n.\end{align}
Going through the arguments of the proof one verifies the error term in \eqref{n260} is locally uniformly bounded.
\end{proof}

\noindent Combining Theorem \ref{specthm} and  Theorem \ref{n3.2ter} one obtains
\begin{cor1}\label{corollary3.2}  (i) For any $q\in H^N_0,\, N\in \mathbb{Z}_{\geq 0}$, 
\begin{align}\label{17a}\tau_n-\mu_n&= \langle q,\cos2\pi nx\rangle +\frac{1}{n^{N+1}}\ell^2_n 
\end{align} 
where the error term is uniformly  bounded on bounded sets of potentials in $H^N_0$.

(ii) For any $N\in \mathbb{Z}_{\geq 0}$, there exists an open neighbourhood $W_N\subseteq  H^N_{0,\mathbb{C}}$ of $H^N_0$ so that \eqref{17a} holds on $W_N$ with a locally uniformly bounded error term.  
\end{cor1}

\section{Appendix A: Infinite products}
In this Appendix we provide asymptotic estimates for infinite products of complex numbers needed to prove the claimed asymptotic estimates of spectral quantities.  For a given sequence $(a_m)_{m\geq 1}$ in $\mathbb{C}$, the infinite product $\prod_{m\geq 1}(1+a_m)$ is said to converge if the sequence $\left(\prod_{1 \leq m\leq M}(1+a_m)\right)_{M\geq 1}$ is convergent. In such a case we set \[\prod_{m\geq 1}(1+a_m)= \lim_{M\rightarrow \infty} \prod_{1 \leq m\leq M}(1+a_m).\] It is said to be absolutely convergent if $\prod_{m\geq 1}(1+|a_m|) < \infty$.
Note that an absolutely convergent infinite product is convergent. A sufficent condition for $\prod_{m\geq 1}(1+a_m)$ being absolutely convergent is that $\|a\|_{\ell^1}:=\sum_{m\geq 1}|a_m| <\infty$. Indeed, as $0\leq \log(1+x) \leq x$ for any $x\geq 0$ one has \begin{equation}\label{def2.1} \prod_{m\geq 1}(1+|a_m|)= \exp\left(\sum_{m\geq 1} \log (1+|a_m|)\right) \leq \exp\left(\sum_{m\geq 1}|a_m|\right).
\end{equation}
We will improve on the estimates of infinite products of Appendix L in \cite{KaP} by using (a version of) the discrete Hilbert transform, defined for an arbitrary sequence $\alpha= (\alpha_m)_{m\geq 1} \in \ell^2_\mathbb{C}$ by
\[ H\alpha:= \left(\sum_{ m\neq n}\alpha_m\left(\frac{1}{n-m}+\frac{1}{n+m}\right)\right)_{n\geq 1}.
\] The following result is due to Hilbert -- see e.g. \cite{HLP}, p 213 for a proof.
\begin{lem1}\label{lemma2.1} H defines a bounded linear operator on $\ell^2_\mathbb{C}$, with  $\|H\|\leq 2\pi$. 
\end{lem1}
\noindent Later we will need the following auxilary result.
\begin{lem1}\label{lemm2.2a} For any $\ell^1$-sequence $(a_m)_{m\geq 1}\subseteq \mathbb{C}$ with $|a_m|\leq \frac{1}{2}$ for any $ m\geq 1$, one has
\[\left|\prod_{ m\geq 1}(1+a_m)-1\right| \leq |A|e^S+|B|e^{S+S^2}
\] where $A=\sum_{m\geq 1} a_m,\, B=\sum_{m\geq 1} |a_m|^2$, and $S=\sum_{m\geq 1} |a_m|$.
\end{lem1}
\begin{proof}As $|a_m|\leq\frac{1}{2}$, the logarithm $\log (1+a_m)$ is well defined and one has
\begin{align*}\prod_{ m\geq 1}(1+a_m)&=\exp\left({\sum_{m\geq 1}\log{(1+ a_{m})}}\right)\\&= e^A\exp\left({\sum_{m\geq 1}\log{(1+ a_{m})}-a_m}\right). \end{align*}
The estimate $\left|\log(1+z)-z\right|\leq \left|z\right|^2$,  $z\in\mathbb{C}$ with $\left|z\right| \leq \frac{1}{2}$, then leads to the following bound for   $R:=\sum_{m\geq 1}(\log(1+ a_{m})-a_m)$,
\[|R|\leq \sum_{m\geq 1}|a_m|^2= B\leq S^2.\]
With  $|A|\leq S$ it then follows that \begin{align*}\left|\prod_{ m\geq 1}(1+a_m)-1\right|= \left|(e^A-1)+e^A(e^R-1)\right| \leq \left|(e^A-1)\right|+e^S\left|(e^R-1)\right|.
\end{align*} As $|e^z-1|\leq |z|e^{|z|}$, $z\in\mathbb{C}$, this leads to the claimed statement, \[\left|\prod_{ m\geq 1}(1+a_m)-1\right|\leq
|A|e^S+|B|e^{S+S^2}.\]
\end{proof}

\noindent To state the results  on infinite products coming up in our study we need to introduce some more notation. 
First let us introduce the notion of isolating neighbourhoods. We say that the discs $U_n\subseteq \mathbb{C},\,n\geq 1$, are a family of \textit{isolating neighbourhoods with parameters} $n_0\geq 1,\, r>0,\,\rho >0$ if the $U_n$'s are mutually disjoint open discs with centers $z_n\in \mathbb{R}$ satisfying $z_1<z_2<... \,,$ so that 
  \begin{align}\label{3bis}U_n\subseteq  D_r^n := \{\lambda \in \mathbb{C} |\; |\lambda- n^2\pi^2|< r\pi^2\}\quad\forall n\geq 1\end{align} and
\[U_n=D_r^n \quad\forall n\geq n_0+1\] so that for any $n,m\geq 1$ 
\begin{align}\label{5bis}|\lambda-\mu|\geq \frac{1}{\rho} |n^2-m^2|\quad \forall \lambda \in U_n,\forall \mu \in U_m.\end{align}
We remark that the results stated below continue to hold for a weaker notion of isolating neighbourhoods, but they suffice for our purposes.
Let $a^0:=(n^2\pi^2)_{n\geq1}$. 
\begin{prop1}\label{corollary2.3}
Assume that $(U_n)_{n\geq 1}$ is a sequence of isolating neighbourhoods with parameters $n_0,r,\rho.$ Then for arbitrary sequences $(a_m)_{m\geq 1},(b_m)_{m\geq 1}\subseteq \mathbb{C}$ with $\alpha := a-a^0,\beta:= b-a^0$ in $\ell^2_\mathbb{C}$, $b_m \in U_m$ for any $m\geq 1$, and any sequence     of complex numbers $\Lambda:=(\lambda_n)_{n\geq 1}$ with $\lambda_n\in U_n\;$ for any $n\geq 1$
\[f_n(\lambda_n)= \prod_{m\neq n}\frac{a_m-\lambda_n}{b_m-\lambda_n}= 1+ \frac{1}{n}\ell^2_n\]  uniformly on bounded subsets of $\alpha,\beta \in \ell^2_\mathbb{C}$ with $b_m\in U_m$ for any $m\geq 1$ and uniformly in $\Lambda$ with $\lambda_n\in U_n$ for any $n\geq 1$. More precisely, \[\sum_{n\geq 1} n^2|f_n(\lambda_n)-1|^2\leq K^2_{\alpha,\beta,\Lambda}\] with a constant $K_{\alpha,\beta,\Lambda}>0$ which  can be chosen uniformly for bounded subsets of $\alpha,\beta \in \ell^2_\mathbb{C}$ with $b_m\in U_m$ for any $m\geq 1$ and uniformly for $\Lambda=(\lambda_n)_{n\geq 1}$ with $\lambda_n \in U_n$ for any $n\geq 1$.
\end{prop1}
\begin{proof}
For any $n,m\geq 1$ with $n\neq m$, \[\frac{a_m-\lambda_n}{b_m-\lambda_n}= 1+ a_{nm}\quad \text{and}\quad a_{nm}= \frac{a_m-b_m}{b_m-\lambda_n}.\]
As by assumption \[|b_m-\lambda_n|\geq \frac{1}{\rho}|m^2-n^2|\quad \forall m\neq n\] and $|m^2-n^2|\geq n|m-n|$, one has by the Cauchy-Schwarz inequality 
\begin{align}\label{5ter}\sum_{m\neq n}|a_{nm}|\leq \rho \sum_{m\neq n}\frac{|a_m-b_m|}{|m^2-n^2|}\leq \frac{\pi \rho}{n}\|a-b\|\end{align} where we used that 
\[\sum_{m\neq n}\frac{1}{|m-n|^2}\leq 2 \sum_{k\geq 1}\frac{1}{k^2}=\frac{\pi^2}{3}.\]
By \eqref{def2.1} this leads to the following estimate
\[\left|\prod_{m\neq n}\frac{a_m-\lambda_n}{b_m-\lambda_n}\right|\leq \prod_{m\neq n}(1+|a_{nm}|) \leq \exp (\frac{\pi \rho}{n}\|a-b\|). \] Now let us consider the asymptotics of $\prod_{m\neq n}\frac{a_m-\lambda_n}{b_m-\lambda_n}$ as $n\rightarrow \infty$. Note that for any $m\neq n$, 
\[\left|\frac{a_m-b_m}{b_m-\lambda_n}\right|\leq \rho \frac{\|a-b\|}{|m^2-n^2|}\leq \frac{\rho}{n}\|a-b\|.\] Choose $n_1\geq n_0$ so that $\frac{\rho}{n_1}\|a-b\|\leq\frac{1}{2}.$ Hence for $n\geq n_1$, Lemma \ref{lemm2.2a} can be applied to $(a_{nm})_{m\neq n}$. To obtain the claimed estimates we need to show that $A_n:= \sum_{m\neq n}a_{nm}= \frac{1}{n}\ell^2_n$ and $B_n:= \sum_{m\neq n} |a_{nm}|^2=\frac{1}{n}\ell^2_n.$ We begin by estimating $A_n$. Note that for $n\geq 1$ and $m\neq n$, 
\begin{align}\label{2.5a} a_{nm} = \frac{b_m-a_m}{n^2\pi^2-m^2\pi^2}+ \frac{b_m-a_m}{n^2\pi^2-m^2\pi^2}\left(\frac{n^2\pi^2-m^2\pi^2}{\lambda_n-b_m}-1\right).
\end{align}
The two terms on the right hand side of the latter identity are treated separately. Use that $(n^2-m^2)^{-1}=\frac{1}{2n}\left(\frac{1}{n-m}+\frac{1}{n+m}\right)$ to conclude from Lemma \ref{lemma2.1} that
\begin{align}\label{2.6a} \left(\sum_{n\geq 1} n^2\left|\sum_{m\neq n}\frac{b_m-a_m}{n^2\pi^2-m^2\pi^2}\right|^2\right)^{\frac{1}{2}} \leq \|b-a\|.
\end{align} To estimate the second term on the right hand side of \eqref{2.5a}, use that $b_m,\lambda_m\in U_m$ so that in view of \eqref{3bis}
\begin{align*}\left|\frac{n^2\pi^2-m^2\pi^2}{\lambda_n-b_m}-1\right| &\leq \frac{|b_m-m^2\pi^2|+|\lambda_n -n^2\pi^2|}{|\lambda_n-b_m|}\\&\leq \rho\frac{2r\pi^2}{|n^2-m^2|},
\end{align*} to obtain the following estimate for the weighted $\ell^1$-norm -- and hence the weighted $\ell^2$-norm -- 
\begin{align}\nonumber &\sum_{n\geq 1}n \sum_{m\neq n} \left|\frac{b_m-a_m}{n^2\pi^2-m^2\pi^2}\right|\left|\frac{n^2\pi^2-m^2\pi^2}{\lambda_n-b_m}-1\right|\\\nonumber \leq & 2 r \rho\sum_{n\geq 1} \sum_{m\neq n}\frac{n}{n^2}\frac{|b_m-a_m|}{(n-m)^2}
\\\nonumber \leq &
2r\rho \left(\sum_{n\geq 1} \sum_{m\neq n}\frac{1}{n^2}\frac{1}{|n-m|^2}\right)^{\frac{1}{2}}\left(\sum_{n\geq 1} \sum_{m\neq n}\frac{|b_m-a_m|^2}{|n-m|^2}\right)^{\frac{1}{2}} \\\label{2.7a} \leq &4\pi r\rho \|b-a\|.
\end{align}
The estimates for the $B_n$'s are simpler as we do not need to split $a_{nm}$. Indeed, as \[B_n \leq \sum_{m\neq n}\rho^2 \frac{|a_m-b_m|^2}{|n^2-m^2|^2}\] we get the following estimate for its weighted $\ell^1$-norm -- and hence for its weighted $\ell^2$-norm --
\begin{eqnarray}\nonumber \sum_{n\geq 1} 
n^2 B_n & \leq &\rho^2 \sum_{n\geq 1} \frac{n^2}{n^2}
\sum_{m\neq n} \frac{|a_m-b_m|^2}{(n-m)^2}\\\nonumber
&\leq  &\rho^2 \|b-a\|^2 
2\sum_{k\geq 1}\frac{1}{k^2} \\\label{2.8a} &\leq & 4\rho^2 \|b-a\|^2. \end{eqnarray} The claimed estimates then follow from Lemma \ref{lemm2.2a},  \eqref{5ter}, and the statement on the uniformity of the estimates follow from the explicit bounds  \eqref{2.6a}--\eqref{2.8a}.
\end{proof}

\begin{cor1}\label{corollary2.4}Assume that $(U_n)_{n\geq 1}$ is a sequence of isolating neighbourhoods with parameters 
$n_0,r,\rho.$ Then for any sequence $a=(a_m)_{m\geq 1}$ with $ a-a^0\in \ell^2_\mathbb{C}$ and any sequence
 $\Lambda =(\lambda_n)_{n\geq 1}$ with $\lambda_n\in U_n$ for any $n\geq 1$, 
the infinite product $\prod_{m\neq n} \frac{a_m-\lambda_n}{m^2\pi^2}$ is absolutely convergent for any  $n\geq 1$ and 
\[\prod_{m\neq n}\frac{a_m-\lambda_n}{m^2\pi^2}=  \frac{(-1)^{n+1}}{2}+ \frac{1}{n}\ell^2_n\] 
uniformly on bounded subsets of $a-a^0 \in\ell^2_\mathbb{C}$ and  uniformly with respect to $\Lambda$ with $\lambda_n\in U_n $ for any $ n\geq 1$.
\end{cor1}
\begin{proof}
For any $m\geq 1, n\geq 1$ write
\[\frac{a_m-\lambda_n}{m^2\pi^2}=1+ \frac{a_m-m^2\pi^2-\lambda_n}{m^2\pi^2}\] and
\[\sum_{m\geq 1} \left|\frac{a_m-m^2\pi^2-\lambda_n}{m^2\pi^2}\right|\leq 
\frac{\|\alpha\|+|\lambda_n|}{\pi^2}\sum_{m\geq 1}\frac{1}{m^2}= \frac{\|\alpha\|+|\lambda_n|}{6},\] where $\alpha:= a-a^0$. Hence according to \eqref{def2.1}, $\prod_{m\neq n}\frac{a_m-\lambda_n}{m^2\pi^2}$ is absolutely convergent and bounded in terms of $\|\alpha\|$ and $|\lambda_n|$. It remains to estimate the product for $n\geq n_0+1$. 
Recall that $\frac{\sin \sqrt\lambda}{\sqrt{\lambda}}$ has the product expansion,
\[\frac{\sin \sqrt\lambda}{\sqrt{\lambda}}=\prod_{m\geq 1}\frac{m^2\pi^2-\lambda}{m^2\pi^2}.\] Hence for any $\lambda_n\in U_n$ with $n\geq n_0+1$ one has
\begin{align}\label{2.9a} \prod_{m\neq n}\frac{a_m-\lambda_n}{m^2\pi^2}=\frac{\sin \sqrt{\lambda_n}}{n^2\pi^2-\lambda_n} \frac{n^2\pi^2}{\sqrt{\lambda_n}}\prod_{m\neq n}\frac{a_m-\lambda_n}{m^2\pi^2-\lambda_n}.\end{align}
In order to apply Proposition \ref{corollary2.3} to the product $\prod_{m\neq n}\frac{a_m-\lambda_n}{m^2\pi^2-\lambda_n}$  we replace for $1\leq m \leq n_0$ the disc $U_m$ by
the disc with center $m^2\pi^2$ and radius $1$. Then the parameters $n_0$ and $r$ can be left as is whereas $\rho$ is replaced by $\rho_1\geq \rho$ so that
\[|m^2\pi^2 -\lambda_n|\geq \frac{1}{\rho_1}|m^2-n^2|\quad\forall \lambda_n \in U_n,\, \forall m\neq n,\, \forall n\geq n_0+1.\]
By Proposition \ref{corollary2.3} \begin{align}\label{2.10a} \prod_{m\neq n}\frac{a_m-\lambda_n}{m^2\pi^2-\lambda_n}= 1 + \frac{1}{n}\ell^2_n 
\end{align}
where the asmptotics are uniform in the sense stated there.  Now let us estimate the remaining terms of the right hand side of \eqref{2.9a}. For this purpose write
\[\frac{\sin \sqrt{\lambda_n}}{n^2\pi^2-\lambda_n} = (-1)^{n+1}\frac{\sin( \sqrt{\lambda_n}-n\pi)}{\sqrt{\lambda_n}-n\pi}\frac{1}{\sqrt{\lambda_n}+n\pi}. \]
Note that $\sqrt{\lambda_n}=n\pi + O(\frac{1}{n})$ and hence by Taylor expansion,
\begin{align*}\sin (\sqrt{\lambda_n}-n\pi) &= (\sqrt{\lambda_n}-n\pi)\left(1+O\left(\frac{1}{n^2}\right)\right)\\ \frac{1}{\sqrt{\lambda_n}+n\pi}&= \frac{1}{2n \pi}\left(1+O\left(\frac{1}{n^2}\right)\right)\\\frac{1}{\sqrt{\lambda_n}}&= \frac{1}{n \pi}\left(1+O\left(\frac{1}{n^2}\right)\right).
\end{align*}
Altogether one obtains in this way
\begin{align}\label{2.11a}\frac{\sin \sqrt{\lambda_n}}{n^2\pi^2-\lambda_n}\frac{n^2\pi^2}{\sqrt{\lambda_n}}= \frac{(-1)^{n+1}}{2}+O\left(\frac{1}{n^2}\right).
\end{align}
One easily sees that the error term $O(\frac{1}{n^2})$ in \eqref{2.11a} can be bounded by $C\frac{1}{n^2}$ where $C>0$ is a constant, only depending on $r$ and $n_0$. Combining estimates \eqref{2.10a} and \eqref{2.11a} then leads to the claimed statement.

\end{proof}

\section{Appendix B: Asymptotics of $\lambda_n$ for real potentials} \label{asymptoticsofperiodicreal}
The aim of this Appendix is to prove the following asymptotics of the periodic eigenvalues of $-d_x^2+q$ for $q$ \textit{real valued}. 
\begin{thm1}\label{reallambdathm} For any $N\in \mathbb{Z}_{\geq 0}$ and $q\in H^N_0$,
\begin{align}\label{n300} \lambda_{2n} =m_n +|\langle q, e^{2\pi i n x} \rangle| + \frac{1}{n^{N+1}}\ell^2_n
\end{align}
and 
\begin{align}\label{n301} \lambda_{2n-1} =m_n -|\langle q, e^{2\pi i n x} \rangle| + \frac{1}{n^{N+1}}\ell^2_n
\end{align}
uniformly on bounded subsets of $H^N_0$. Here $m_n$  is the expression given by \eqref{12a}.
\end{thm1}
\begin{rem1}
For $N\geq 1$, the  asymptotics \eqref{n300}-\eqref{n301} are due to Marchenko  \cite{Ma}, whereas for the uniformity statement we could not find any reference.
\end{rem1}
\begin{proof}We again use the special solutions $z_N(x,\nu)$ introduced in Section \ref{sectionspectialsolutions}. Following Marchenko we represent these special solutions for $|\nu|$ sufficiently large by an exponential function , $z_N(x,\nu)= \exp \left(i\nu x +\int_0^x \sigma(t,\nu) dt\right) $, where 
\begin{align}\label{d131bis}\sigma(t,\nu)=\sum_1^N\frac{s_k(t)}{(2i\nu)^k} + \frac{\sigma_N(t,\nu)}{(2i\nu )^N}\end{align}	
and $s_k(t)$ are given as in Section \ref{sectionspectialsolutions}. Recall that \[z_N(x,\nu)= \exp \left( i\nu x+\sum_{k=1}^N\int_0^x\frac{s_k(t)}{(2i\nu)^k} dt\right) + \frac{r_N(x,\nu)	}{(2\pi i\nu)^{N+1}}.\] Hence for
$|\nu|$ large,
\[\int^x_0 \frac{\sigma_N(t,\nu)}{(2i\nu )^N}d t=\log\left(1+\exp \left( -i\nu x-\sum_{k=1}^N\int_0^x\frac{s_k(t)}{(2i\nu)^k} dt \right) \cdot \frac{r_N(x,\nu)}{(2\pi i\nu)^{N+1}}\right) .\]
Furthermore, as $r_N(0,\nu)=0$ and $r'(0,\nu)=0$ it follows that 
\[\sigma_N(0,\nu)=0.\]
As $z'_N(x,\nu)= z_N(x,\nu)(i\nu +\sigma(x,\nu))$,  the determinant of the solution matrix $Y_N(x,\nu)$, given by
\[Y_N(x,\nu)= \begin{pmatrix} z_N(x,-\nu) & z_N(x,\nu)\\z'_N(x,-\nu) &z'_N(x,\nu) 
\end{pmatrix}\]  can be computed
to be \[\det Y_N(x,\nu)= z_N(x,-\nu)z_N(x,\nu)\cdot w(x,\nu)\] where
\begin{align}\label{d131ter} w(x,\nu)= 2i\nu +\sigma(x,\nu)- \sigma(x,-\nu )  .\end{align}
Note that $\det Y_N(0,\nu)= w(0,\nu)$. As the Wronskian is $x$- independent we get 
\[ z_N(1,-\nu)z_N(1,\nu)\cdot w(1,\nu)=w(0,\nu)\]
or \begin{align}\label{n303}z_N(1,-\nu)z_N(1,\nu) =\frac{w(0,\nu)}{w(1,\nu)}.
\end{align} This identity allows to express $z_N(1,-\nu)$ in terms of $z_N(1,\nu)$,
\begin{align}
\label{n304}z_N(1,-\nu) =\frac{w(0,\nu)}{w(1,\nu)}\frac{1}{z_N(1,\nu)}.
\end{align}
Further note that the fundamental matrix is given by $Y_N(x,\nu)Y_N(0,\nu)^{-1}$. The condition that $\nu^2$ be a periodic eigenvalue can be expressed by 
\[\chi_p(\nu)=\det \left(Y_N(1,\nu)Y_N(0,\nu)^{-1}-Id_{2\times 2} \right)=0\]
whereas anti-periodic eigenvalues are characterized by \[\chi_{ap}(\nu)=\det \left(Y_N(1,\nu)Y_N(0,\nu)^{-1}+Id_{2\times 2}\right)=0.\]
The two cases are treated similary and so we concentrate on the first one.
Clearly, $\det \left(Y_N(1,\nu)Y_N(0,\nu)^{-1}-Id_{2\times 2}\right)= \frac{\det (Y_N(1,0)-Y_N(0,\nu))}{\det (Y_N(0,\nu))}$ and
\begin{align*} & \det (Y_N(1,0)-Y_N(0,\nu))\\ =& \det \begin{pmatrix} z_N(1,-\nu)-1 & z_N(1,\nu)-1 \\
z'_N(1,-\nu)-z_N'(0,-\nu) & z'_N(1,\nu)-z_N'(0,\nu)
\end{pmatrix}\\
= & \Big(z_N(1,-\nu)-1\Big)\Big( z_N(1,\nu)\left(i\nu +\sigma(1,\nu)\Big)- \left( i\nu+ \sigma(0,\nu)\right)\right)
\\ & -\Big(z_N(1,\nu)-1\Big)\Big( z_N(1,-\nu)(-i\nu +\sigma(1,-\nu))- ( -i\nu+ \sigma(0,-\nu))\Big)\\
 =& -\Big(2i\nu+ \sigma(1,\nu)-\sigma(0,-\nu)\Big)z_N(1,\nu) +\Big(-2i\nu+ \sigma(1,-\nu)-\sigma(0,\nu)\Big)z_N(1,-\nu)\\
&+ z_N(1,-\nu)z_N(1,\nu)\Big(2i\nu +\sigma(1,\nu)-\sigma(1,-\nu)\Big) + 2i\nu +\sigma(0,\nu) -\sigma(0,-\nu). 
\end{align*} 
Introduce $G(\nu):= 2i\nu +\sigma(1,\nu )-\sigma(0,-\nu)$. With $z_N(1,-\nu)z_N(1,\nu) =\frac{w(0,\nu)}{w(1,\nu)}$ one then gets
\begin{align*} \det (Y_N(1,0)-Y_N(0,\nu))=& -G(\nu)z_N(1,\nu)+ G(-\nu)z_N(1,-\nu)\\
&+\frac{w(0,\nu)}{w(1,\nu)}\cdot w(1,\nu) +w(0,\nu).
\end{align*} 
The equation $\chi_p(\nu)=0$ is thus equivalent to 
\[G(\nu) z_N(1,\nu)-2w(0,\nu)-G(-\nu)z_N(1,-\nu)=0.\]
Dividing this equation by $w(1,\nu)z_N(1,-\nu)=\frac{w(0,\nu)}{z_N(1,\nu)}$ leads to the following  quadratic equation 
\[\frac{G(\nu)}{w(0,\nu)}z_N(1,\nu)^2-2z_N(1,\nu)-\frac{G(-\nu)}{w(1,\nu)}=0\]
or \begin{align}\label{n310}z_N(1,\nu)=\frac{w(0,\nu)}{G(\nu)}\pm \frac{w(0,\nu)}{G(\nu)}\sqrt{1+\frac{G(\nu)G(-\nu)}{w(0,\nu)w(1,\nu)}}.
\end{align}
We now bring the right hand side of this identity in a more convenient form. For this purpose introduce 
\begin{align}\label{d134bis}D(\nu)=\frac{\sigma(1,\nu)-\sigma(0,\nu)}{w(0,\nu)}.\end{align}
Using that $w(0,\nu)=2i\nu +\sigma(0,\nu)-\sigma(0,-\nu)$ we get 
\begin{align*}G(\nu) = 2i\nu + \sigma(1,\nu)- \sigma (0,-\nu) 
= w(0,\nu)+\sigma(1,\nu) -\sigma(0,\nu)
\end{align*} or 
\begin{align}\label{n311} G(\nu)= w(0,\nu)(1+D(\nu)).\end{align}
As $w(0,-\nu)=-w(0,\nu)$ by the definition of $w(0,\nu)$ one gets 
\begin{align*} \frac{G(\nu)G(-\nu)}{w(0,\nu)w(1,\nu)}=& \frac{w(0,\nu)(1+D(\nu))w(0,-\nu)(1+D(-\nu))}{w(0,\nu)w(1,\nu)}\\
=& - \frac{1+D(\nu)+D(-\nu)+D(\nu)D(-\nu)}{1+\frac{w(1,\nu)-w(0,\nu)}{w(0,\nu)}}.
\end{align*}
As \begin{align*} w(1,\nu)-w(0,\nu)=& 2i\nu +\sigma(1,\nu)-\sigma(1,-\nu) -2i\nu + \sigma(0,\nu)+\sigma(0,-\nu)\\
=& w(0,\nu)(D(\nu)+D(-\nu))
\end{align*} 
 we then get \begin{align*}
\frac{G(\nu)G(-\nu)}{w(0,\nu)w(1,\nu)}=& -\frac{1+D(\nu)+D(-\nu)+ D(\nu)D(-\nu)}{1+D(\nu)+D(-\nu)}
\\ =& -1-\frac{D(\nu)D(-\nu)}{1+D(\nu)+D(-\nu)}
\end{align*}
Substituting this identity as well as \eqref{n311} into \eqref{n310} yields
\begin{align}\label{n312} z_N(1,\nu)=\frac{1}{1+D(\nu)}	\left( 1\pm i\sqrt{\frac{D(\nu)D(-\nu)}{1+D(\nu)+D(-\nu)}}\right).
\end{align}
We have to take a closer look at $D(\nu)$. Recall that $D(\nu)=  \frac{\sigma(1,\nu)-\sigma(0,\nu)}{w(0,\nu)}$ where \[\sigma(x,\nu)= i\nu + \sum_1^N\frac{s_k(x)}{(2i\nu)^k} + \frac{\sigma_N(x,\nu)}{(2i\nu)^N}.\] As the $s_k$'s are 1-periodic and $\sigma_N(0,\nu)=0$ it follows that \[\sigma(1,\nu)-\sigma(0,\nu)=\frac{\sigma_N(1,\nu)}{(2i\nu)^N}.\]
Writing 
\[w(0,\nu) =2i\nu + \sigma(0,\nu)-\sigma(0,-\nu)=2i\nu \left(1+ \frac{\sigma(0,\nu)-\sigma(0,-\nu)}{2i\nu} \right)\] then leads to
\begin{align*}D(\nu)=& \frac{\sigma_N(1,\nu)}{(2i\nu)^{N+1}}\left(1+\frac{\sigma(0,\nu)-\sigma(0,-\nu)}{2i\nu}\right)^{-1}\\
=& \frac{\sigma_N(1,\nu)}{(2i\nu)^{N+1}}\left(1+O\left(\frac{\sigma(0,\nu)-\sigma(0,-\nu)}{2i\nu}\right)\right).
\end{align*} 
As $\sigma(0,\nu)= \sum_1^N\frac{s_k(0)}{(2i\nu)^k}$ one has \[\sigma(0,-\nu)-\sigma(0,\nu)= \sum_{k=1}^N\frac{s_k(0)}{(2i\nu)^k}((-1)^k-1)=O\left(\frac{1}{\nu}\right)\]
and hence 
\begin{align}\label{n315} D(\nu)= \frac{\sigma_N(1,\nu)}{(2i\nu)^{N+1}}\left(1+O\left(\frac{1}{\nu^2}\right)\right). \end{align}
In view of this estimate for $D(\nu)$ we can take the logarithm of \eqref{n312} for $\nu=\nu_n$ where $\nu_n=\sqrt[+]{\lambda_{2n}}$ or $\nu_n= \sqrt[+]{\lambda_{2n-1}}$ with $n$ \textit{even}. (Recall that the periodic eigenvalues of $-d_x^2+q$ when considered on $[0,1]$ are given by $\lambda_0 < \lambda_3 \leq \lambda_4<\lambda_7\leq \lambda_8< \dots$ .) Due to the asymptotics $\lambda_{2n},\lambda_{2n-1}=n^2\pi^2 +\ell^2_n$ it follows that $\nu_n= n \pi+\frac{1}{n}\ell^2_n$. Taking the logarithm of 
\[z_N(1,\nu_n)= \exp \left(i(\nu_n-n\pi)+ \sum_1^N \int_0^1 \frac{s_k(t)}{(2i\nu_n)^k}dt + \int_0^1 \frac{\sigma_N(t,\nu_n)}{(2i\nu_n)^N} dt\right).\]
As $n$ is even, the identity \eqref{n312} together with Lemma \ref{lemA.3}  leads to 
\begin{align*} &i(\nu_n-n\pi) -i\sum_{1\leq 2l+1\leq N}\frac{(-1)^la_{2l+1}}{(2\nu_n)^{2l+1}}+ \int_0^1 \frac{\sigma_N(t,\nu_n)}{(2i\nu_n)^N} dt
\\=& -\log \left(1+D(\nu_n)\right) + \log \left(1\pm i\frac{\sqrt{D(\nu_n)D(-\nu_n)}}{\sqrt{1+D(\nu_n)+D(-\nu_n)}}\right)
\\=& -\frac{\sigma_N(1,\nu_n)}{(2i\nu_n)^{N+1}} +\frac{1}{2} \frac{\sigma(1,\nu_n)^2}{(2i\nu_n)^{2N+2}} +O\left(\left(\frac{1}{n^{N+3}}\right)\right)	  \pm 
i\frac{\sqrt{D(\nu_n)D(-\nu_n)}}{\sqrt{1+D(\nu_n)+D(-\nu_n)}}\\ &+ \frac{1}{2} \frac{{D(\nu_n)D(-\nu_n)}}{{1+D(\nu_n)+D(-\nu_n)}}
+O\left(\left(\frac{{D(\nu_n)D(-\nu_n)}}{{1+D(\nu_n)+D(-\nu_n)}}\right)^{3/2}\right).
\end{align*}
By the estimate \eqref{n315} and the expansion $(1+x)^{-1/2}= 1-x/2+3x^2/8+O(x^3)$ one gets 
\begin{align*}
\left(1+D(\nu_n) +D(-\nu_n)\right)^{-1/2}=& 1-\frac{1}{2} \frac{1}{(2i\nu_n)^{N+1}}\Big(\sigma_N(1,\nu_n)
\\& + (-1)^{N+1}\sigma_N(1,-\nu_n)\Big)+ O\left(\frac{1}{n^{2N+2}}\right)\end{align*}
As $\sqrt{D(\nu_n)D(-\nu_n)}= O\left(\frac{1}{n^{N+1}}\right)$ it then follows that 
\begin{align*}\frac{\sqrt{D(\nu_n)D(-\nu_n)}}{\sqrt{1+D(\nu_n)+D(-\nu_n)}}= {\sqrt{D(\nu_n)D(-\nu_n)}} 
-\frac{1}{2} \frac{\sqrt{D(\nu_n)D(-\nu_n)}}{(2i\nu_n)^{N+1}}\\ \cdot\Big(\sigma_N(1,\nu_n)+ (-1)^{N+1}\sigma_N(1,-\nu_n)\Big) + O\left(\frac{1}{n^{N+3}}\right).
\end{align*}
Combining the estimates above leads to 
\begin{align}\nonumber  &\nu_n -n\pi -\sum_{1\leq 2l+1\leq N}\frac{(-1)^la_{2l+1}}{(2\nu_n)^{2l+1}} -i\int_0^1 \frac{\sigma_N(t,\nu_n)}{(2i\nu_n)^N} dt 
\\\nonumber =& \; 	i \frac{\sigma_N(1,\nu_n)}{(2i\nu_n)^{N+1}} -\frac{i}{2}\frac{\sigma_N(1,\nu_n)^2}{(2i\nu_n)^{2N+2}} 
 \pm \Big({\sqrt{D(\nu_n)D(-\nu_n)}}\\\nonumber  &- \frac{1}{2} \frac{\sqrt{D(\nu_n)D(-\nu_n)}}{(2i\nu_n)^{N+1}} \left(\sigma_N(1,\nu_n)+(-1)^{N+1}\sigma_N(1,-\nu_n)\right)\Big) 
\\\label{n318} &-\frac{i}{2} D(\nu_n)D(-\nu_n) +O\left(\frac{1}{n^{N+3}}\right).
\end{align}
We need now to consider $\sigma_N(1,\nu_n)$ and $\int_0^1 \sigma_N(t,\nu_n)dt$.
For  this purpose introduce the following notation
\begin{align*} A:=&\exp\left(i\nu-i\sum_{1\leq 2l+1\leq N}\frac{(-1)^la_{2l+1}}{(2\nu )^{2l+1}}\right)\\ 
B:=& \exp\left(\int_0^1 \frac{\sigma_N(t,\nu)}{(2i\nu)^N} dt\right).
\end{align*}
By the definition of $z_N$ and $\sigma_N$, \[z_N(1,\nu)=A\cdot B= A+\frac{r_N(1,\nu)}{(2i\nu)^{N+1}}.\]
Hence \[B= 1+A^{-1}\frac{r_N(1,\nu)}{(2i\nu)^{N+1}}\]
and, by taking logarithm, 
\begin{align}\label{n3200}\int_0^1\frac{\sigma_N(t,\nu)}{(2i\nu)^N}dt= A^{-1}\frac{r_N(1,\nu)}{(2i\nu)^{N+1}}-\frac{1}{2}A^{-2}\frac{r_N(1,\nu)^2}{(2i\nu)^{2N+2}}
+O\left(\frac{1}{\nu^{3N+3}}\right).
\end{align}
Furthermore, with $z_N(1,\nu)=A\cdot B$, 
\begin{align*}z'_N(1,\nu) =& z_N(1,\nu)\left(i\nu + \sum_1^N \frac{s_k(0)}{(2i\nu)^k}+\frac{\sigma_N(1,\nu)}{(2i\nu)^{N}}\right)\\ 
=& A\left(i\nu +\sum_1^N \frac{s_k(0)}{(2i\nu)^k}\right) +\frac{r'_N(1,\nu)}{(2i\nu)^{N+1}}
\end{align*}
or, as $z_N(1,\nu)-A= r_N(1,\nu)/(2i\nu)^{N+1}$,
\[z_N(1,\nu)\frac{\sigma_N(1,\nu)}{(2i\nu)^N}= -\frac{r_N(1,\nu)}{(2i\nu)^{N+1}}\left(i\nu +\sum_1^N \frac{s_k(0)}{(2i\nu)^k}\right) +\frac{r'_N(1,\nu)}{(2i\nu)^{N+1}}.\]
It leads to the formula \begin{align}\label{n320} z_N(1,\nu)\sigma_N(1,\nu) =
-\frac{r_N(1,\nu)}{2}- \frac{r_N(1,\nu)}{2i\nu}\sum_1^N \frac{s_k(0)}{(2i\nu)^k}+ \frac{r'_N(1,\nu)}{2i\nu}.
\end{align}  
Recall that by Proposition \ref{lemA.2} (keep in mind that n is \textit{even})
\[r_N(1,\pm\nu_n)= a_{N+1} -(\pm 2in\pi)^N\langle q,e^{\mp 2in\pi x} \rangle \pm \frac{1}{2in\pi} a_{N+2} +\frac{1}{n}\ell^2_n\]
and 
\[r'_N(1,\pm\nu_n)= \pm in\pi a_{N+1} \pm in\pi(\pm 2in\pi)^N\langle q,e^{\mp 2in\pi x} \rangle + \frac{ a_{N+2}}{2} +\frac{1}{n}\ell^2_n.\]
Using that $a_1=0$ and $\nu_n-n\pi =\frac{1}{n} \ell^2_n$ it follows that 
\[z_N(1,\nu_n)=\exp \left(i(\nu_n -n\pi)+ \sum_{1\leq	2l+1 \leq N} \frac{a_{2l+1}}{(2i\nu_n)^{2l+1}}\right)= 1+\frac{1}{n}\ell^2_n.\]
Hence \begin{align}\nonumber z_N(1,\pm\nu_n)\sigma_N(1,\pm \nu_n) = &-\frac{a_{N+1}}{2} +\frac{(\pm 2in\pi)^N}{2}\langle q,e^{\mp 2in\pi x} \rangle \mp  \frac{1}{2in\pi} \frac{a_{N+2}}{2} +\frac{1}{n}\ell^2_n \\\nonumber &+\frac{a_{N+1}}{2} +\frac{(\pm 2in\pi)^N}{2}\langle q,e^{\mp 2in\pi x} \rangle \pm  \frac{1}{2in\pi} \frac{a_{N+2}}{2} +\frac{1}{n}\ell^2_n \end{align}
and thus
\begin{align}	\label{n322} 
 z_N(1,\pm\nu_n)\sigma_N(1,\pm \nu_n)= {(\pm 2in\pi)^N}\langle q,e^{\mp 2in\pi x} \rangle +\frac{1}{n}\ell^2_n .
\end{align}  

\noindent As $z_N(1,\pm\nu_n)=1+O\left(\frac{1}{n}\right)$ 
\begin{align}\label{mai79bis} \frac{\sigma_N(1,\pm\nu_n)}{(\pm 2i\nu_n)^{N+1}}= \frac{1}{\pm 2in\pi} \langle q,e^{\mp 2\pi i nx} \rangle +\frac{1}{n^{N+2}}\ell^2_n. \end{align}
Hence also \[\frac{\sigma_N(1,\pm \nu_n)^2}{(\pm 2i\nu_n)^{2N+2}}= \frac{1}{n^{2N+2}}\ell^2_n\]
and  
\[\frac{\sqrt{D(\nu_n)D(-\nu_n)}}{(2i\nu_n)^{N+1}}\left( \sigma_N(1,\nu_n) +(-1)^{N+1}\sigma_N(1,-\nu_n)\right)= \frac{1}{n^{2N+2}}\ell^2_n\]
as well as 
\[D(\nu_n)D(-\nu_n)= \frac{\sigma_N(1,\nu_n)\sigma_N(1,-\nu_n)}{(2\nu_n)^{2N+2}}\left(1+O\left(\frac{1}{n^2}\right)\right)=\frac{1}{n^{2N+2}}\ell^2_n.\]
In view of all this and \eqref{n3200}, \eqref{n318} reads  
\begin{align*} \nu_n = & n\pi + \sum_{1\leq 2l+1\leq N}\frac{(-1)^la_{2l+1}}{(2\nu_n)^{2l+1}} + iA^{-1} \frac{r_N(1,\nu_n)}{(2i\nu_n)^{N+1}}  \\
& -\frac{i}{2} A^{-2} \frac{r_N(1,\nu_n)^2}{(2i\nu_n)^{2N+2}} 
+ \frac{1}{2n\pi} \int_0^1 q(x) e^{2\pi inx} dx 
\\ & +\sqrt{ D(\nu_n)D(-\nu_n))} +\frac{1}{n^{N+2}}\ell^2_n .
\end{align*}
But \[A^{-1}=1+ O\left(\frac{1}{n^2}\right)\] hence
\[iA^{-1} \frac{r_N(1,\nu_n)}{(2i\nu_n)^{N+1}} = i\frac{a_{N+1}}{(2i\nu_n)^{N+1}}- \int_0^1 \frac{q(x)e^{2\pi i nx}}{2n\pi}dx + i\frac{a_{N+2}}{(2i\nu_n)^{N+2}}\]
and using that $a_{N+1}=0$ in the case $N=0$ we get for any $N\geq 0$ that
\[ A^{-2}\frac{r_N(1,\nu_n)}{(2i\nu_n)^{2N+2}} = O\left(\frac{1}{n^{N+3}}\right) .\] We then get in view of Lemma \ref{lemA.3}
\begin{align*} \nu_n = & n\pi + \sum_{1\leq 2l+1\leq N+2}\frac{(-1)^la_{2l+1}}{(2\nu_n)^{2l+1}}  
\pm \sqrt{ D(\nu_n)D(-\nu_n))} +\frac{1}{n^{N+2}}\ell^2_n .
\end{align*}

Note that the periodic spectrum of $-d_x^2+q$ is real and hence $\sqrt[+]{\lambda_{2n}},\sqrt[+]{\lambda_{2n-1}}$ are real for $n$ sufficiently large. It then follows from Lemma \ref{dlem12.1} below that
\[\sqrt[+]{D(\nu_n)D(-\nu_n)}= |D(\nu_n)|\]
But by \eqref{n315},  $D(\nu_n)= \frac{\sigma_N(1,\nu_n)}{(2in\pi)^{N+1}}+O\left(\frac{1}{n^{N+3}}\right)$ and by \eqref{mai79bis}
\[ \frac{\sigma_N(1,\nu_n)}{(2in\pi)^{N+1}}= \frac{1}{2in\pi}\int_0^1 q(x) e^{2in\pi x} dx + \frac{1}{n^{N+2}}\ell^2_n.\]
Hence, for $q$ real we have by Cauchy's inequality  
\begin{align*}& \left||\frac{\sigma_N(1,\nu_n)}{(2in\pi)^{N+1}}|- \frac{1}{2n\pi} |\int_0^1 q(x) e^{2\pi in x} dx|\right|
\\ \leq & \left|\frac{\sigma_N(1,\nu_n)}{(2in\pi)^{N+1}}-\frac{1}{2in\pi} \int_0^1 q(x) e^{2\pi in x} dx\right| = \frac{1}{n^{N+2}}\ell^2_n. 
\end{align*}
We therefore have established that 
\begin{align*} \nu_n = & n\pi + \sum_{1\leq 2l+1\leq N+2}\frac{(-1)^la_{2l+1}}{(2\nu_n)^{2l+1}} \pm \left|\frac{1}{2n\pi} \int_0^1 q(x) e^{2\pi in x} dx\right| +\frac{1}{n^{N+2}}\ell^2_n .
\end{align*}
Arguing as for the asymptotics of the Dirichlet eigenvalues one gets
\begin{align*} \nu_n = & n\pi + \sum_{1\leq 2l+1\leq N+2}\frac{(-1)^lb_{2l+1}}{(2n\pi)^{2l+1}} \pm \left|\frac{1}{2n\pi} \int_0^1 q(x) e^{2\pi in x} dx\right| +\frac{1}{n^{N+2}}\ell^2_n.
\end{align*}
Then for $q\in H^N_0$ 
\[\nu_n^2= m_n \pm \left|\int_0^1 q(x) e^{2\pi in x} dx\right| +\frac{1}{n^{N+1}}\ell^2_n \]
where $m_n$ is given by \eqref{12a}. 
Going through the arguments of the proof one concludes that the error term has the claimed uniformity property.  Finally we need to determine the signs $\pm$ in the asymptotics of $\lambda_{2n}$ and $\lambda_{2n-1}$. It is convenient to introduce $\lambda^+_n=\lambda_{2n},\; \lambda_n^-=\lambda_{2n-1}.$ Choose $\epsilon_n^{\pm} \in \{-1,1\}$ so that 
\[\lambda_n^{\pm} =m_n+ \epsilon_n^{\pm}\left|\int_0^1 q(x) e^{2\pi in x} dx\right|+\frac{1}{n^{N+1}}\ell^2_n.\]
We note that $\epsilon_n^{\pm}$ are not uniquely determined and that as $\lambda_n^-\leq  \lambda_n^+$, we may choose 
$\epsilon_n^{\pm}$ so that $\epsilon_n^-\leq \epsilon_n^+$. 
We claim that $\epsilon_n^{\pm}$ can be chosen in such a way that $\epsilon_n^{+}=1$ and $\epsilon_n^{-}=-1$ for any $ n\geq 1$. Indeed, let $J=\{n\geq 1|\; \epsilon_n^{+}=\epsilon_n^{-}\}$. If $J$ is finite we can change $\epsilon_n^{+},\epsilon_n^{-}$ as claimed without changing the asymptotics. If $J$ is infinite, then the $n$'s in $J$ form a subsequence $(n_k)_{k\geq 1}$ in $\mathbb{N}$ such that $\epsilon_{n_k}^{+}=\epsilon_{n_k}^{-}$ for any $ k\geq 1$. As 
\[\tau_{n_k}= \frac{\lambda_{n_k}^+ +\lambda_{n_k}^-}{2} = m_{n_k} +\epsilon_{n_k}^{+}\left|\int_0^1 q(x) e^{2\pi in_k x} dx\right|+ \frac{1}{n_k^{N+1}} \ell^2_k\]
it then follows  from the asymptotics of the $\tau_n$'s that \[\left|\int_0^1 q(x) e^{2\pi in_k x} dx\right|_{k\geq 1}= \frac{1}{n_k^{N+1}}\ell^2_k\]
and hence we can again change the $\epsilon_{n_k}^{\pm}$ so that now $\epsilon_{n_k}^{\pm}=\pm 1$.
This proves the claimed asymptotics.
\end{proof}

\begin{lem1}\label{dlem12.1} For any $\nu \in \mathbb{R}$ and any $q\in L^2_0$ (in particular $q$ real valued)
\begin{enumerate} \item [(i)]$z_N(x,-\nu)=\overline{z_N(x,\nu)}\quad \forall x\in \mathbb{R}$.
\item[(ii)] $\sigma(x,-\nu)=\overline{\sigma(x,\nu)}\quad \forall x\in \mathbb{R}$
\item[(iii)] $D(-\nu) = \overline{ D(\nu)} $.
\end{enumerate} 
\end{lem1}
\begin{proof}(i) 
It is easy to check that $z(x,-\nu)$ and $\overline{z(x,\nu)}$ both satisfy the equation 
\begin{align}\label{star} -y'' +qy =\nu^2 y.\end{align}
Further recall that 
\[z(0,-\nu) =1=\overline{z(0,\nu)}\] and 
\[z'(0,-\nu) =-i\nu +\sum_1^N\frac{s_k(0)}{(-2i\nu)^k}= \overline{i\nu +  \sum_1^N\frac{s_k(0)}{(2i\nu)^k}}=\overline{z'(0,\nu)} \]
By the uniqueness of solutions of \eqref{star} with given initial values it then follows that \[z(x,-\nu)=\overline{z(x,\nu)}\quad \forall x\in \mathbb{R},\forall \nu \in \mathbb{R}.\] 

\noindent (ii) By  \eqref{84bis}-\eqref{A.2} one sees that $s_k(x)$ is real for any $x\in \mathbb{R}$ and any $1\leq k\leq N.$ Hence by (i) and the definition \eqref{d131bis} of $\sigma(t,\nu)$ it follows that
\[\sigma(t,-\nu)= \overline{\sigma(t,\nu)}.\]

\noindent (iii) In view of the definition of $D(\nu)$ (\eqref{d134bis},\,\eqref{d131ter}), the claimed identity follows from (ii).
\end{proof}


\begin{thebibliography}{123}
%\bibitem{BG} D. Bambusi, B. Gr\'ebert: {\em Birkhoff normal form for PDE's with tame modulus},
%Duke Math. J. 135 (2005), 507-567.

\bibitem{FM} H. Flaschka, D. McLaughlin: {\em Canonically conjugate variables for the Korteweg-de Vries equation and Toda lattices with periodic boundary conditions}, Progress Theor.Phys. 55(1976), 438-456.

%\bibitem{GKP} B. Gr\'{e}bert, T. Kappeler, J. P\"oschel: {\em Normal form theory of the NLS equation}, preliminary version, arxiv.org/abs/0907.3938v1 2009



\bibitem{HLP} G. Hardy, J. Littlewood, G. P\'olya: {\em Inequalities}, Cambridge University Press, paperback edition, reprinted 1994.



%\bibitem{Ka} T. Kappeler: {\em Fibration of the phase-space for the Korteweg-de Vries equation}, 
 % Ann. Inst. Fourier 41 (1991), 539-575.

%\bibitem{KM} T. Kappeler, M. Makarov: {\em On action-angle variables for the second Poisson bracket}, Comm. in Math. Phys., 214, 3 (2000), 651-677.

%\bibitem{KaMa} T. Kappeler, M. Makarov: {\em On Birkhoff 
%coordinates for KdV}, Ann. H. Poincar\'e 2 (2001), 807 - 856.

%\bibitem{KaMi} T. Kappeler, B.Mityagin: {\em Estimates for periodic and Dirichlet eigenvalues of the Schr\"odinger operator}, SIAM J. Math Anal. 33 (2001), p 113-152. 

%\bibitem{KMT} T. Kappeler, C. M\"ohr, P. Topalov: {\em Birkhoff coordinates for KdV on phase spaces of distributions}, 
 %Selecta Math. (N.S.) 11, no 1(2005), 37-98. 
 
   \bibitem{KaP} T. Kappeler, J. P\"oschel: {\em KdV \& KAM}, 
  Ergebnisse Math. u. Grenzgebiete, Springer, Berlin, 2003.

%\bibitem{KP1} T. Kappeler, J. P\"oschel: {\em On the periodic KdV equation in weighted Sobolev spaces}, 
%  Ann. H. Poincar\'e  (C) Non Linear Anal., 26, no 3 (2009), 841-853.



%\bibitem{KST1} T. Kappeler, B. Schaad, P. Topalov: {\em mKdV and its Birkhoff coordinates}, Physica D: Nonlinear Phenomena, 237, 10-12 (2008), 1655-1662.

\bibitem{KST3} T. Kappeler, B. Schaad, P. Topalov: {\em Qualitative features of periodic solutions of KdV}, preprint.


%\bibitem{KST2} T. Kappeler, B. Schaad, P. Topalov: {\em On the asymtotics of the Birkhoff map of NLS}, %in preparation.

%\bibitem{KST} T. Kappeler, F. Serier, P. Topalov: {\em On the symplectic phase space of KdV}, 
 % Proc. Amer. Math. Soc. 136 (2008), 1691-1698.

%\bibitem{KT1} T. Kappeler, P. Topalov: {\em Global fold structure of the Miura map on $L^2(\mathbb{T},\mathbb{R})$}, 
 % Int. Math. Research Notices (2004), 2039-2068.


%\bibitem {Ku} S. Kuksin: {\em Damped driven KdV and effective equation for long-time behaviour of its solution}, preliminary version, 2009.

%\bibitem {KP} S. Kuksin, G. Perelman: {\em Vey theorem in infinite dimensions
 %   and its application to KdV}, Discrete and Continuous Dynamical Systems  Ser. A, 27, no 1, (2010),  1-24.

%\bibitem {KuPi} S. Kuksin, A. Piatnitski: {\em Khasminskii - Whitham averaging for randomly perturbated KdV equation}, J. Math. Pures Appl. 89 (2008), 400-428.

%\bibitem {Lev} B. Levitan: {\em Inverse Sturm-Liouville problems}, VSI. Zeist, 1987.


\bibitem {Ma} V. Marchenko: {\em Sturm-Liouville operators and applications}, Birkh\"auser, Basel, 1986.
   
%\bibitem{MKV} H. McKean, E. Vaninsky: {\em Action-angle variables for the cubic Schroedinger equation}, Comm. Pure Appl. Math. 50 (1997), 489-562.

%\bibitem{Poschel} J. P\"oschel: {\em Hill's potentials in weighted Sobolev spaces and their spectral gaps}, to appear in Math. Annalen. DOI 10.1007/s 00208-010-0513-7.
  
 \bibitem{PT} J. P\"oschel, E. Trubowitz: {\em Inverse spectral theory}, Academic Press, Boston, 1987.
 
\bibitem{SS} A.  Savchuk, A. Shkalikov: {\em On the eigenvalues of the Sturm-Liouville operator with potentials from Sobolev spaces}, Math. Notes 80, no 6 (2006), 864-884.

%\bibitem{VN} A. Veselov, S. Novikov: {\em Poisson brackets and complex tori}, Proc. Steklov Institue 1985, issue 3,  53-65.
  

 

   



\end{thebibliography}
\end{document}